\documentclass[12pt]{amsart}
\usepackage{amsmath}
\usepackage{amsfonts}
\usepackage{amssymb}
\usepackage[mathscr]{eucal}
\usepackage{graphicx}
\makeatletter
\def\numberwithin#1#2{\@ifundefined{c@#1}{\@nocnterrr}{%
  \@ifundefined{c@#2}{\@nocnterr}{%
  \@addtoreset{#1}{#2}%
  \toks@\expandafter\expandafter\expandafter{\csname the#1\endcsname}%
  \expandafter\xdef\csname the#1\endcsname
    {\expandafter\noexpand\csname the#2\endcsname
     .\the\toks@}}}}
\makeatother
\numberwithin{equation}{section}

\newtheorem{theorem}{Theorem}
\numberwithin{theorem}{section}

\newtheorem{corollary}[theorem]{Corollary}

\newtheorem{definition}[theorem]{Definition}
\newtheorem{example}[theorem]{Example}
\newtheorem{lemma}[theorem]{Lemma}

\newtheorem{proposition}[theorem]{Proposition}
\newtheorem{remark}[theorem]{Remark}

\newenvironment{rmk}{\begin{remark} \em}{\end{remark}}

\newcommand{\tr}{\operatorname{tr}}
\newcommand{\Rc}{\operatorname{Rc}}

\newcommand{\SO}{{\rm SO}}
\newcommand{\SU}{{\rm SU}}
\newcommand{\Sp}{{\rm Sp}}
\newcommand{\Spin}{{\rm Spin}}

\newcommand{\so}{\mbox{${\mathfrak s \mathfrak o}$}}

\newcommand{\m}{\mbox{${\mathfrak m}$}}

\newcommand{\C}{\mbox{${\mathbb C}$}}

\newcommand{\HH}{\mbox{${\mathbb H}$}}

\newcommand{\PP}{\mbox{${\mathbb P}$}}

\newcommand{\R}{\mbox{${\mathbb R}$}}

\begin{document}
\title{Ancient Solutions on Bundles with Non-abelian Structural Group}

\author{Peng Lu}
\address{Department of Mathematics, University of Oregon, Eugene, OR 97403-1222}
\email{penglu@uoregon.edu}
\thanks{}

\author{Y. K. Wang}
\address{Department of Mathematics and Statistics, McMaster
     University, Hamilton, Ontario, L8S 4K1, CANADA}
\email{wang@mcmaster.ca}
\thanks{P.L. is partially supported by Simons Foundation through Collaboration Grant 229727.
Y.K.W. is partially supported by NSERC Grant No. OPG0009421.}

\date{\small revised \today}

\begin{abstract}
We generalize the ancient solutions of the Ricci flow on certain principal ${\rm SO}(3)$ bundles
over compact quaternionic K\"ahler manifolds constructed by Bakas, Kong, and Ni  to
certain $\R\PP^3$ fibre bundles over a product of two compact quaternionic K\"ahler manifolds.
The ancient solutions are of Type I, $\kappa$-noncollapsed, and have positive Ricci curvature.
Analogous solutions on even-dimensional non-K\"ahlerian bundle spaces are also constructed.

\smallskip


\end{abstract}

\maketitle

\section{\bf Introduction}

Ancient solutions for the Ricci flow arise naturally when one blows up the finite time
singularities of Ricci flow solutions on closed manifolds. The work in \cite{BKN12} and
\cite{DHS12} shows that ancient solutions are interesting objects in geometric analysis.

In \cite{LW16} we generalized the examples in \cite[\S 5]{BKN12} of ancient solutions on principal
circle bundles over Fano K\"ahler-Einstein manifolds to the situation of principal torus bundles
over a finite product of Fano K\"ahler-Einstein manifolds. In this paper we shall generalize
in a similar fashion the ancient solutions given by Theorem 6.2 in \cite{BKN12} on certain
principal ${\rm SO}(3)$ bundles over compact quaternionic K\"ahler manifolds (abbreviated as
QK manifolds in the sequel).

To describe our results, let us choose QK manifolds $(M_i^{4n_i},g_i), \, i =1, \cdots, m$,
with $m \geq 2$. Recall that for each such space there is an associated principal
$\operatorname{SO}(3)$ bundle $P_i$ over $M_i$ whose associated $\R^3$ bundle is the analog
of the anti-canonical bundle for the K\"ahler-Einstein case. Taking the product of these principal bundles
we can form an associated $\frac{\SO(3) \times \cdots \times \SO(3)}{\Delta \SO(3)}$ {\it fibre bundle}
$\bar{P}$ over the product  $M_1 \times \cdots\times M_m$ (see \S \ref{subsec 2.1 intro to QK} for further
details). By \cite{Wa92}, there exist connection type Einstein metrics on $\bar{P}$. As in
\cite{BKN12}, we exploit the fact that the Ricci flow preserves the family of connection type
metrics, and so the search for ancient solutions of connection type on $\bar{P}$ is reduced to the search for
special flow lines for a reduced system of ordinary differential equations ($\operatorname{ODE}$s in short).

In this paper we will construct continuous families of $\kappa$-noncollapsed
type I ancient solutions when $m=2$ and for a special class of bundles when $m\geq 3$.
Below we give a description of the $m=2$ case, and refer the reader to
Theorem \ref{thm sec 5 m> 3} for  our results in the  $m\geq 3$ case.
Note that $\frac{\SO(3) \times \SO(3)}{\Delta \SO(3)}$ is diffeomorphic to
$\R\PP^3$, and so the ancient solutions exist on $\R\PP^3$ fibre bundles.

\vskip .3cm
\noindent{\bf Main Theorem.} {\em  There are $($i$)$ a continuous $1$-parameter
family of ancient solutions and $($ii$)$ two distinct
ancient solutions of the Ricci flow on the $\R\PP^3$ fibre bundle $\bar{P}$  over an
arbitrary product $M_1 \times M_2$ of QK manifolds with positive scalar curvature.
These solutions are of type I when time $t \rightarrow - \infty$ and are
$\kappa$-noncollapsed at all scales for some $\kappa > 0$. Their rescaled backwards
limits are Einstein metrics. Furthermore, the solutions in both $($i$)$ and  $($ii$)$
have positive Ricci curvature for $t$ close to $-\infty$. }

\vskip .2cm

Technically, the difference between our work and that in \cite{BKN12} is that the $\operatorname{ODE}$
system in \cite{BKN12} is two-dimensional, has a first integral, and hence is integrable.
For our system, there is no longer an obvious first integral, and the two Einstein metrics on
the bundle $\bar{P}$ are only established implicitly.
To get ancient solutions whose rescaled backwards limits are these Einstein metrics,
by  the Hartman-Grobman Theorem, it suffices for us to show the existence of unstable directions
of the linearization of our system about these Einstein metrics.  Thus one needs
to achieve this using appropriate estimates for their location since they are not given
explicitly. Furthermore, for some of these ancient solutions,
we will consider their singularity formulation in finite time, which requires us
to find suitable bounded sets and differential inequalities and use them to show that
the solutions stay in these sets.  The geometric properties of the solutions are then
established from various estimates derived from the $\operatorname{ODE}$s.

Compared to the work in \cite{LW16}, the analysis here is more delicate because the bundle
$\bar{P}$ admits two Einstein metrics. One of the Einstein metrics takes the place of the
origin in the system considered in \cite{LW16}, which is a source for the Ricci flow.
The origin is then an additional stationary solution for the system in this paper. The fact
that the origin is actually a sink for the Ricci flow means that we do not get collapsed solutions
anymore. This can perhaps be regarded as an indication through the Ricci flow that
bundles with non-abelian structural groups are more rigid than those with abelian structural
groups.

The ancient solutions in \S \ref{subsec 5 QK m geq 3 same structure}
 include continuous families which exist on even dimensional
non-K\"ahlerian manifolds. The only examples of this type in the literature that we are aware of
are a $2$-parameter family of left-invariant solutions found in \cite{La13} on each of
the compact connected simple Lie groups except $\Sp(n)$.

As in \cite{LW16}, we also obtain without additional work the existence of pseudo-Riemannian
ancient and immortal solutions if we replace some or all of the QK factors of positive scalar
curvature by ones with negative scalar curvature (see Remark \ref{pseudo-Rflow}).

The following is an outline of the rest of this paper.
In \S \ref{sec 2, ODE and linearizaion} we derive the $\operatorname{ODE}$ for the Ricci flow and
establish properties about the linearization of the  $\operatorname{ODE}$ at Einstein metrics
(Proposition \ref{prop fixed point asyp type}).
In \S \ref{sec 3 ancient sol existence metric behav}  we prove the existence of ancient
solutions on the  $\R\PP^3$ fibre bundles over a product of two QK manifolds (Theorem
\ref{thm psi b one two}), and discuss their metric behavior as time
$t \rightarrow -\infty$ and their $\kappa$-noncollapsed property  (Theorem \ref{thm growth ai bi QK Omega 1}
and \ref{thm growth ai bi QK Omega 3}, and Corollary \ref{cor noncollapsing}).
In \S \ref{subsec 4 QK curva prop} we prove a general result (Theorem \ref{thm 1 no choice}) 
which enables us to use the metric tensor limits of the previous section to conclude 
(without any computations) that  the corresponding ancient solutions are of Type I as time approaches
$-\infty$. We also show that the ancient solutions for the $m=2$ case have positive Ricci curvature
when time is close to $-\infty$ (Theorem \ref{thm anc sol curva prop}).
In \S \ref{subsec 5 QK m geq 3 same structure} we prove the existence of ancient
solutions on a special subfamily of fibre bundles over a product of an arbitrary number of QK manifolds
(Theorem \ref{thm sec 5 m> 3}) and discuss their properties.

\section{\bf Ricci flow on fibre bundles over products of two QK manifolds
and its linearization at Einstein metrics}\label{sec 2, ODE and linearizaion}

In this section we derive a system of ODEs which is equivalent to the Ricci flow
equation for a family of connection type metrics on certain
$\frac{\operatorname{SO}(3) \times \cdots \times \operatorname{SO}(3)) }{ \Delta \operatorname{SO}(3)}$
fibre bundles over a product of $m$ arbitrary QK manifolds. We then specialize to the
$m=2$ case, and determine the stationary solutions of the $\operatorname{ODE}$ system
and  establish properties of the linearization of
the  $\operatorname{ODE}$ system at these solutions which we will need
in the later sections.

\subsection{Some facts about QK manifolds} \label{subsec 2.1 intro to QK}

We first mention some standard facts about QK manifolds, referring the reader to
\cite{Be87} and \cite{Sa89} for details. Let $(M^{4n},g), n\geq 2,$ be a QK manifold, i.e.,
its holonomy group is contained in the subgroup $\Sp(n)\cdot \Sp(1)$ of $\SO(4n)$. ($\Sp(n) \cdot \Sp(1)$
denotes more precisely the group $(\operatorname{Sp}(n) \times \operatorname{Sp}(1) )/\Delta\mathbb{Z}_2$,
where $\Delta\mathbb{Z}_2$ denotes the diagonally embedded central subgroup.)
We further assume that its holonomy group does not lie in $\Sp(n)$. Then it must be
de Rham irreducible and Einstein with nonzero scalar curvature.
A QK manifold with positive or negative scalar curvature will be called a
{\em positive} or {\em negative} QK manifold, respectively. When $n=1$, it is
customary to regard self-dual Einstein manifolds as the $4$-dimensional analogs of QK
manifolds. Those with positive scalar curvature are then known to be isometric to $S^4$
or $\C\PP^2$ by a celebrated theorem of N. Hitchin \cite{Hi81}. We shall include the $n=1$
case in all of our discussions below.

To fix notation for this paper, let $\Rc_g =\Lambda g$ and let $\hat{P}$ denote the
principal $\Sp(n) \cdot \Sp(1)$ holonomy bundle associated to the tangent bundle of $M$.
The Levi-Civita connection of $g$ induces a connection on $\hat{P}$, which upon projection
induces a connection on the principal
\[
((\operatorname{Sp}(n) \times \operatorname{Sp}(1) )/\Delta\mathbb{Z}_2)/\operatorname{Sp}(n)
 \simeq \operatorname{SO}(3)
\]
bundle $P \doteqdot {\hat P}/\Sp(n)$ over $M$. This connection is in fact Yang-Mills,
and has parallel curvature form \cite{MS88}. It is the analogue of the connection
induced by the Levi-Civita connection of a K\"{a}hler-Einstein
manifold on the associated circle bundle of its anti-canonical line bundle.
The $ \operatorname{SO}(3)$ bundle $P$ is in general not spin, i.e., it does not in general admit a lift to an
$\SU(2)$ bundle. Indeed, when $M$ is compact, it is spin precisely when $M$ is a
quaternionic projective space.

Let $(M_i^{4n_i}, g_i), \, i =1,\cdots, m$, be arbitrary QK manifolds with
$\Rc_{g_i} =\Lambda_i\, g_i$. (Note that there is actually no need to assume that $\Lambda_i > 0$
at this stage.) Let $P_i \rightarrow M_i$ be the principal
$\operatorname{SO}(3)$ bundle described above that is associated to $(M_i,g_i)$.
Let $ \Delta \operatorname{SO}(3)$ be the diagonal subgroup of the product group
$\operatorname{SO}(3) \times \cdots \times \operatorname{SO}(3)$.
 We define an
 \[
 \frac{\operatorname{SO}(3) \times \cdots \times \operatorname{SO}(3)}{ \Delta
\operatorname{SO}(3)}
\]
fibre bundle $\bar{P}$ of dimension $\left(\sum_{i=1}^m 4n_i\right) +3(m-1)$ over $\bar{M}
 \doteqdot M_1 \times \cdots \times M_m$ by
\[
\bar{P} \doteqdot (P_1 \times \cdots \times P_m) \times_{\operatorname{SO}(3) \times \cdots \times
\operatorname{SO}(3)} \left ( \frac{\operatorname{SO}(3) \times \cdots \times \operatorname{SO}(3)}{ \Delta
\operatorname{SO}(3)} \right ) \rightarrow \bar{M}.
\]
By the above discussion, each $P_i$  has a Yang-Mills connection whose curvature form is parallel.
The product of these connections gives a connection on the
$\operatorname{SO}(3) \times \cdots \times
\operatorname{SO}(3)$  principal bundle $P_1 \times \cdots \times P_m$.

Using this connection on $P_1 \times \cdots \times P_m$ and choosing vectors
$\vec{a}=(a_1, \cdots, a_m)$ and $\vec{b}=(b_1, \cdots, b_m)$ with {\it all positive}
components, we can construct a family of Riemannian metrics $\bar{g}_{\vec{a},\vec{b}}$
 on $\bar{P}$ uniquely determined by $\vec{a}$ and $\vec{b}$ satisfying the following
conditions:

\medskip

\noindent (B1) $(\bar{P}, \bar{g}_{\vec{a},\vec{b}}) \rightarrow (\bar{M}, g_{\vec{b}})$ is
a Riemannian submersion with totally geodesic fibers, where
$g_{\vec{b}} \doteqdot \sum_{i=1}^m b_i \pi_i^* g_i$ and $\pi_i:
\bar{M} \rightarrow M_i$ is the natural projection map.

\vskip .2cm
\noindent (B2) The restriction of $\bar{g}_{\vec{a},\vec{b}}$ to the typical fiber in $\bar{P}$
is the normal homogeneous metric on $\frac{\operatorname{SO}(3) \times \cdots \times \operatorname{SO}(3)}
{ \Delta \operatorname{SO}(3)} $ induced from the bi-invariant metric
  $a_1 B_{\so(3)} \oplus \cdots \oplus a_m B_{\so(3)}$
on $\operatorname{SO}(3) \times \cdots \times \operatorname{SO}(3)$
where $B_{\so(3)}(U, V) = -\tr(UV)$ for $U, V \in \so(3)$, the Lie algebra of $\SO(3)$.

\medskip

Note that the isotropy representation of the fibre consists of $m-1$ copies of
the adjoint representation of $\SO(3)$, so the space of invariant metrics on
$ \frac{\operatorname{SO}(3) \times \cdots \times
\operatorname{SO}(3)} { \Delta \operatorname{SO}(3)} $ is given by
\[
\operatorname{GL}^+(m-1) / \operatorname{SO}(m-1) \equiv S_+^2(\mathbb{R}^{m-1})
\]
with dimension $\frac{m(m-1)}{2}$. The invariant metrics on $\frac{\operatorname{SO}(3) \times
 \cdots \times  \operatorname{SO}(3)} { \Delta \operatorname{SO}(3)} $ induced from
bi-invariant metrics on  $\operatorname{SO}(3) \times \cdots \times \operatorname{SO}(3)$ as
described in condition (B2) above is an $m$-dimensional subfamily except in the case $m=2$ (see
Remark \ref{effparam}).

\begin{rmk}\label{rem 5.1 QK conj}
It should be mentioned that there is a conjecture, widely believed to be true, that a positive
QK manifold must be isometric to a quaternionic symmetric space of compact type, also known as
a Wolf space. In fact the $4$- and $8$-dimensional cases have been verified in \cite{Hi81} and
\cite{PS91}, respectively. A further rigidity result, valid for all dimensions,
 under the assumption of a non-vanishing second Betti number was obtained in \cite{LS94}.

If the fundamental conjecture about positive QK manifolds is true, then each $M_i$ will
have the form $G_i/(H_i \cdot \Sp(1))$, where $G_i$ is a compact simple Lie group, and so
$P_i = G_i/H_i$. It follows that $\bar{P}$ and the corresponding Ricci flows are homogeneous.
The status of the conjecture aside, we believe that the bundle formalism used in this paper
may be more helpful for understanding ancient solutions beyond the homogeneous category.
Furthermore, as we shall  observe in Remark \ref{pseudo-Rflow}, our analysis also leads
to examples of ancient and immortal solutions of pseudo-Riemannian Ricci flows when some of the
QK factors have negative  scalar curvature.
\end{rmk}

\subsection{Ricci flow on bundle $\bar{P}$:  an $\operatorname{ODE}$ system}
\label{subsec 5.2 metric anciet QK}

Since the connection used to define our metrics $\bar{g}_{\vec{a},\vec{b}}$ is Yang-Mills,
it follows that the only nonzero components of the Ricci tensor of $\bar{g}_{\vec{a},\vec{b}}$
come exclusively from the base or from the fibre, i.e., there are no cross terms involving
both the fibre and the base. To derive the Ricci flow equations for the metrics $\bar{g}_{\vec{a},\vec{b}}$
on $\bar{P}$, we use the computations in \cite{Wa92}, particularly equations (2.1) and (2.2) there.

Consider a family of such Riemannian metrics  $\bar{g}(\tau) \doteqdot \bar{g}_{\vec{a}(\tau), \vec{b}(\tau)}$,
where $\tau \in \mathbb{R}$. First we consider tangential directions along the base.
Let  $\{ X_{k}^{(i)} \}_{k=1}^{4n_i }$ be an orthonormal frame in compact QK manifold $(M_i
^{4n_i}, g_i)$ and let
$\{ \tilde{X}_{k}^{(i)} \}_{k=1}^{4n_i }$ be its horizontal lift to $\bar{P}$.
Using \cite[(2.1)]{Wa92}, it follows immediately that the part of the backwards
Ricci flow equation involving the base directions is given by the equations

\begin{align}
\frac{db_i}{d \tau} &=\frac{d \bar{g}}{d \tau} (  \tilde{X}_{k}^{(i)},  \tilde{X}_{k}^{(i)})
= 2 \Rc_{\bar{g}} (  \tilde{X}_{k}^{(i)},  \tilde{X}_{k}^{(i)})  \notag  \\
& =2\Lambda_i - \frac{6 \Lambda_i^2}{(n_i +2)^2} \cdot \frac{a_i}{b_i} \left (1 -\frac{a_i}{\hat{a}}
\right ), \quad i =1, \cdots, m,  \label{eq quaternion b}
\end{align}
where $1\leq k\leq 4n_i$ is arbitrary but fixed and $\hat{a} \doteqdot \sum_{j=1}^m a_j$.

Next we consider the nonzero components of the Ricci tensor of  $\bar{g}(\tau)$
coming from the fibre directions.  Let $U$ be an element of Lie algebra $\so(3)$ .
 For $i > j$ we define $\tilde{U}_{ji} = (0,\cdots,0,-a_i U, 0, \cdots, 0,a_jU, 0, \cdots, 0)$
 in which $-a_i U$ is at the $j$-th entry and $a_jU$ is at the $i$-th entry.
 For $i < j$ we define $\tilde{U}_{ji} = (0,\cdots,0,a_j U, 0, \cdots, 0, -a_iU, 0, \cdots, 0)$
 in which $a_j U$ is at the $i$-th entry and $-a_iU$ is at the $j$-th entry.
  Then
\[
\tilde{U}_{ji} \in \m_{\vec{a}} \doteqdot \{ (U_1, \cdots, U_m) \in \so(3) \times
\cdots \times \so(3), \,\, \sum_{i=1}^m a_i U_i =0 \},
\]
where  $\m_{\vec{a}}$ is the orthogonal complement of diagonal $\Delta \so(3)$ in
$ \so(3) \times \cdots \times \so(3)$ with respect to the bi-invariant
metric $\oplus_{i=1}^m a_i B_{\so(3)}$. We divide the following discussion into two parts:
$m \geq 3$ and $m=2$.

\vskip .2cm
When  $m \geq 3$, using \cite[(2.2)]{Wa92} we have that for distinct indices $i,  j, k$
\begin{align}
& \frac{d \bar{g}}{d \tau} (\tilde{U}_{ki}, \tilde{U}_{kj} )  =\frac{d a_k}{d \tau}\, a_i a_j B(U,U),
 \label{eq g tilde U} \\
& \Rc_{\bar{g}}(\tilde{U}_{ki},\tilde{U}_{kj}) =\left (\frac{1}{4} + \frac{a_k}{2 \hat{a}}
+ \frac{2n_k \Lambda_k^2 }{(n_k+2)^2} \cdot \frac{a_k^2}{b_k^2} \right ) a_i a_j B(U,U).
\label{eq Rc tilde U}
\end{align}
From the backwards Ricci flow equation $\frac{d \bar{g}}{d \tau} (\tilde{U}_{ki}, \tilde{U}_{kj} )   =
2\Rc_{\bar{g}}(\tilde{U}_{ki}, \tilde{U}_{kj}  )$  we have
\begin{equation}
 \frac{d  a_k }{d \tau} = \frac{1}{2} + \frac{a_k}{ \hat{a}}
+ \frac{4n_k \Lambda_k^2 }{(n_k+2)^2} \cdot \frac{a_k^2}{b_k^2} . \label{eq quaternion a evol}
\end{equation}
It follows that $\hat{a}$ satisfies the equation
\begin{equation}
\frac{d \hat{a}}{d \tau} =\frac{m}{2}+1 + \sum_{k=1}^m 4n_k q_k^2 \cdot
\frac{a_k^2}{b_k^2}, \label{eq d a hat d tau}
\end{equation}
where $q_k \doteqdot \frac{\Lambda_k}{n_k+2}$.

\vskip .2cm
When $m=2$, the complement  $\m_{\vec{a}}$ equals $\{\tilde{U}
\doteqdot (-a_2 U,a_1U): \,  U \in \so(3) \}$.
By a computation similar to that in the $m \geq 3$ case we get the following formulas
\begin{align}
\frac{d \bar{g}}{d \tau} (\tilde{U},  \tilde{U}) =& \left( a_2^2   \frac{d a_1}{d \tau} +a_1^2
 \frac{d a_2}{d  \tau} \right)\, B(U,U), \notag \\
 \Rc_{\bar{g}}(\tilde{U},\tilde{U}) = & \left (\frac{1}{4} + \frac{a_1}{2 (a_1+a_2)}
+ 2n_1q_1^2  \cdot \frac{a_1^2}{b_1^2} \right ) a_2^2 \,B(U,U)    \label{eq Rc tilde U m=2} \\
&  + \left (\frac{1}{4} + \frac{a_2}{2 (a_1+a_2)}+ 2n_2 q_2^2  \cdot
\frac{a_2^2}{b_2^2} \right ) a_1^2 \,B(U,U). \notag
\end{align}
Let $\psi \doteqdot \frac{a_1a_2}{a_1 +a_2}$. It follows from the backwards Ricci flow equation
$\frac{d \bar{g}}{d \tau} (\tilde{U},  \tilde{U}) = 2 \Rc_{\bar{g}}(\tilde{U},\tilde{U})$ and
the relation $\frac{d \psi}{d\tau} = \frac{\psi^2}{a_1^2 a_2^2}   ( a_2^2   \frac{d a_1}{d \tau}
+a_1^2  \frac{d a_2}{d  \tau})$  that
\begin{equation}
 \frac{d \psi}{d \tau} = \frac{1}{2} +4 n_1 q_1^2 \cdot \frac{\psi^2}{b_1^2}
+ 4 n_2 q_2^2  \cdot \frac{\psi^2}{b_2^2}.
  \label{eq m=2 a evol m=2}
\end{equation}

On the other hand, if we set $a_1 = a_2 = 2\psi > 0$ in the expressions of
$\frac{d \bar{g}}{d \tau} (\tilde{U},  \tilde{U})$ and $\Rc_{\bar{g}}(\tilde{U},\tilde{U})$
above, then the fibre part of the  backwards Ricci flow equation becomes (\ref{eq m=2 a evol m=2}).
There is no loss of generality to make this choice since the fibres are the irreducible
symmetric space $\R\PP^3$ when $m=2$, and there is only one effective parameter for the space
of invariant metrics.

\begin{rmk} \label{effparam}
In order to identify the effective parameters for the fibre metrics, one has to break the
symmetry. Notice that the complement $\m_{\vec{a}}$ to the diagonal subalgebra $\Delta \so(3)$ we chose
depended on the biinvariant metric chosen. When $m=2$, it is natural to fix the complement
to be the anti-diagonal $\m \doteqdot \{(U, -U): U \in \so(3) \}$, which is the orthogonal complement when using
$B_{\so(3)} \oplus B_{\so(3)}$ as the background metric.
One checks easily that
$a_1B_{\so(3)} \oplus a_2B_{\so(3)}$ and $a_1^{\prime} B_{\so(3)} \oplus a_2^{\prime} B_{\so(3)}$
induce the same metric on  $\m$
if and only if $a_1 + a_2 = a_1^{\prime} + a_2^{\prime}$.

Alternatively, we can also include $\so(3)$ into the first summand of $\so(3) \oplus \so(3)$
and project the image onto subspace $\m_{\vec{a}}$ using $a_1B_{\so(3)} \oplus a_2B_{\so(3)}$.
Then unit vector $U$ in $(\so(3), B_{\so(3)})$
is mapped to $\frac{-1}{a_1 + a_2}(-a_2 U, a_1 U)$, whose length squared equals to the parameter
$\psi = a_1a_2/(a_1 + a_2)$. This is what we used before in deriving our equations.
\end{rmk}

\subsection{Properties of Einstein metrics on bundle $\bar{P}$ when $m=2$}  \label{subsec 2.3 linearlization}

Note that when $m=2$, $\bar{P}$ is an $\R\PP^3 \approx \frac{\operatorname{SO}(3) \times
\operatorname{SO}(3)} { \Delta \operatorname{SO}(3)} $ fibre bundle over $M_1 \times M_2$.
Let $\bar{g}(\tau) =  \bar{g}_{\vec{a}(\tau), \vec{b}(\tau)}$ be a solution of the backwards
Ricci flow on $\bar{P}$ given in \S \ref{subsec 5.2 metric anciet QK}.
It follows from (\ref{eq m=2 a evol m=2})  and (\ref{eq quaternion b}) that
$\psi(\tau) = \frac{a_1a_2}{a_1 +a_2}$,  $b_1 (\tau)$ and $b_2 (\tau)$ satisfy the following system:
\begin{subequations}
\begin{align}
&  \frac{d \psi}{d \tau} = \frac{1}{2} + 4 n_1 q_1^2 \cdot \frac{\psi^2}{b_1^2}
+ 4 n_2 q_2^2 \cdot \frac{\psi^2}{b_2^2},
  \label{eq m=2 a evol m=2 vers 2} \\
& \frac{db_i}{d \tau}=2(n_i+2)q_i - 6 q_i^2 \cdot  \frac{\psi}{b_i},
 \quad i =1, \, 2,  \label{eq m=2 b 1,2}
\end{align}
\end{subequations}
where $q_i =\frac{\Lambda_i}{n_i+2}$. By the discussion in Remark \ref{effparam},
these equations represent the backwards Ricci flow on $\bar{P}$ for the family of metrics in which
$a_1(\tau)=a_2(\tau) =2\psi(\tau)$.

\begin{rmk}
(i) The difference between the above equations and equations (3.2a)
and (3.2b) in \cite{LW16} for circle bundles is the term $\frac{1}{2}$ in
(\ref{eq m=2 a evol m=2 vers 2}). This term reflects that fact that the
fibers $ \R\PP^3$ now have constant positive curvature.

(ii) The difference between the  above equations and equations (7.2) and (7.3) in \cite{BKN12}
is that the base manifold $\bar{M}$ is reducible with two de Rham factors.
Note that only the product
of two Einstein manifolds with the same Einstein constants is Einstein. So the hypotheses
in Proposition 7.1 in \cite{BKN12} are seldom satisfied in the present situation.
\end{rmk}

To analyse the above system, as in \cite[\S 3.1]{LW16}, we define new dependent variables
$Y_i \doteqdot \frac{\psi}{b_i}, \, i=1, 2$, and a new independent variable
\begin{equation}
u=u(\tau) \doteqdot \int_0^{\tau} \frac{d \zeta }{\psi(\zeta)} .  \label{eq def u via psi}
\end{equation}
Using (\ref{eq m=2 a evol m=2 vers 2}) and (\ref{eq m=2 b 1,2}) we compute that
\begin{align*}
\frac{d Y_i}{du} & =\psi \frac{dY_i}{d \tau} = Y_i \frac{d\psi}{d \tau} - Y_i^2 \frac{d b_i}{d \tau} \\
& = Y_i \left (\frac{1}{2} -  2(n_i+2)q_i Y_i +6q_i^2Y_i^2+4 n_1 q_1^2Y_1^2+ 4n_2 q_2^2Y_2^2
\right  ).
\end{align*}
Let $Y$ denote the vector  $(Y_1,Y_2)$ and let
\begin{align*}
& E(Y) \doteqdot \frac{1}{2} + 4 n_1 q_1^2Y_1^2+ 4n_2 q_2^2Y_2^2, \\
& F_i(Y) \doteqdot  2(n_i +2)q_iY_i -6q_i^2Y_i^2 -E(Y), \quad i=1, 2.
\end{align*}
Equations (\ref{eq m=2 a evol m=2 vers 2}) and (\ref{eq m=2 b 1,2}) become
\begin{subequations}
\begin{align}
& \frac{1}{\psi}  \frac{d \psi}{du} =E(Y),  \label{eq m=2 a evol m=2 vers u} \\
&  \frac{d Y_i}{du} = - Y_i F_i(Y), \quad i=1, 2.  \label{eq m=2 b 1,2 vers u}
\end{align}
\end{subequations}

Given a solution $Y(u)$ of (\ref{eq m=2 b 1,2 vers u}), $\psi (u)$ is determined by
(\ref{eq m=2 a evol m=2 vers u}) up to a multiplicative constant. By the relation
 $d \tau = \psi(\tau(u))du$,  $\tau$ is recovered up to the same multiplicative
 constant since in (\ref{eq def u via psi}) we arranged for $u=0$ to correspond to
 $\tau = 0$. Then $b_i(u)$ can be found from $\frac{\psi}{Y_i}$, and  we get
$(\psi(\tau),b_1(\tau),b_2(\tau))$ and the corresponding backwards flow solution $\bar{g}(\tau)$.
Notice that a multiplicative constant actually amounts to a parabolic rescaling
of $\bar{g}(\tau)$. Hence different integration constants do not lead to new solutions
of the backwards Ricci flow. We have thus proved that each solution $Y(u)$ gives
rise to a single solution $\bar{g}(\tau)$ of the backwards Ricci flow.
Therefore we will focus the solutions of (\ref{eq m=2 b 1,2 vers u}) from now on.

\smallskip

As a preliminary step we have the following lemma about the stationary solutions, two
of which correspond to Einstein metrics on $\bar{P}$.
For the rest of this section  we will \textbf{assume} that $n_1 \leq n_2$
(by interchanging $M_1$ and $M_2$ if necessary).

\begin{lemma} \label{lem station sol QK}
The zeros of vector field  $ (-Y_1 F_1(Y), -Y_2 F_2(Y))$  are

\noindent{$($i$)$}  the origin,

\noindent{$($ii$)$} the points $v_1= \left(\frac{1}{ (4n_1+6) q_1},0 \right), \, \tilde{v}_1
=\left(\frac{1}{2 q_1},0\right),
\, v_2= \left(0, \frac{1}{ (4n_2+6) q_2}\right), \, \tilde{v}_2
=\left(0, \frac{1}{2 q_2}\right)$,

\noindent{$($iii$)$} the Einstein points $\xi=(\xi_1,\xi_2)$ and $\eta =(\eta_1,\eta_2)$.
We further have $ 0< \eta_1 < \xi_1$ and  $0< \eta_2 < \xi_2$.
\end{lemma}

\begin{proof}
$($i$)$  If both $Y_1$ and $Y_2$ are zero, we get the origin.

\vskip .1cm
 $($ii$)$ This occurs when exactly one of $Y_1$ or $Y_2$ is zero. If $Y_2 =0$,
then $F_1(Y_1,0)=0$.  We observe that equation
\[
F_1(Y) =-(4n_1+6)q_1^2Y_1^2 +2(n_1+2)q_1Y_1 -4n_2q_2^2Y_2^2 -\frac{1}{2} =0
\]
represents an ellipse in $Y$-space with center at $\left(\frac{n_1+2}{(4n_1 +6) q_1}, 0\right)$.
The ellipse intersects the $Y_1$-axis at $v_1 =\left(\frac{1}{(4n_1 + 6)q_1}, 0\right)$
 and $\tilde{v}_1 =\left(\frac{1}{2q_1}, 0\right)$, i.e., solutions of the equation $F_1(Y_1,0) =0$.
So for points lying in the ellipse we have $ \frac{1}{(4n_1 + 6)q_1} \leq Y_1 \leq \frac{1}{2q_1}$
and $|Y_2| \leq \frac{n_1 + 1}{2q_2 \sqrt{n_2(4n_1 + 6)}}$.

Likewise, if $Y_1 =0$, then $F_2(0, Y_2) = 0$. $F_2(Y) = 0$ is an ellipse with center at
$\left(0, \frac{n_2+2}{(4n_2 +6) q_2}\right)$ and whose intersection with the
 $Y_2$-axis occur at $v_2 =\left(0, \frac{1}{(4n_2 + 6)q_2}\right)$ and
 $\tilde{v}_2 =\left(0, \frac{1}{2q_2}\right)$. So for points lying in the ellipse
 we have $ \frac{1}{(4n_2 + 6)q_2} \leq Y_2 \leq \frac{1}{2q_2}$ and
$|Y_1| \leq \frac{n_2 + 1}{2q_1\sqrt{n_1(4n_2 + 6)}}$.

\vskip .1cm
$($iii$)$ If both $Y_1$ and $Y_2$ are non-zero, we get $F_1(Y)=F_2(Y) =0$.
Now we consider the relation between the solutions of equations $F_1(Y)=F_2(Y)=0$ and the
Einstein metrics on $\bar{P}$ of the form $\bar{g}_{\vec{a},\vec{b}}$. We have
\begin{equation}
2 (n_1+2)q_1 Y_1 - 6 q_1^2  Y_1^2 =E(Y) = 2(n_2+2)q_2 Y_2 - 6 q_2^2  Y_2^2.
\label{eq E y1 y2 a}
\end{equation}
 Set  $E(Y)\doteqdot 2 \Lambda \psi$ for some $\Lambda >0$.
Hence for $i=1,2,$
\begin{equation}
\frac{1}{4} + \frac{2n_1 \Lambda_1^2}{(n_1+2)^2} Y_1^2 + \frac{2n_2 \Lambda_2^2}
 {(n_2+2)^2}Y_2^2 = \Lambda \psi = \Lambda_i Y_i - \frac{ 3 \Lambda_i^2}{(n_i+2)^2}  Y_i^2 .
 \label{eq E y1 y2 b}
\end{equation}
So $\psi$ and $b_i$ satisfy the equations for an Einstein metric
(\cite[(2.3)]{Wa92}) and we have $\Lambda = \frac{\Lambda_i}{b_i} -\frac{3 \Lambda_i^2}{(n_i+2)^2}
\cdot \frac{\psi}{b_i^2}$.
It is known from \cite[p.313]{Wa92} that there are exactly two solutions
$Y=\xi=(\xi_1,\xi_2)$ and  $Y=\eta=(\eta_1,\eta_2)$ of the above equations which represent
Einstein metrics on $\bar{P}$.

\vskip .1cm
To deduce the inequalities in (iii), let
 \[
 y_0 \doteqdot y_0(Y)   \doteqdot \frac{1}{2E(Y)}, \quad y_i  \doteqdot y_i(Y)  \doteqdot
 4 (n_i+2) q_i y_0(Y) Y_i, \quad i=1, 2.
 \]
 Then the equations $F_1(Y)= F_2(Y) =0$ together with the definition of $E(Y)$ give rise to
an equivalent system of equations for the variables $y_0, y_1, y_2$,
\begin{subequations}
\begin{align}
& 3(y_0-1) +2n_1 (y_1-1) + 2n_2(y_2-1) =0   , \label{eq einstein y0} \\
&     y_0 = \frac{3}{4(n_i+2)^2} \cdot \frac{y_i^2}{y_i -1}, \quad i=1, 2. \label{eq einstein y i}
\end{align}
\end{subequations}
Actually (\ref{eq einstein y i}) follows from (\ref{eq E y1 y2 a}) and (\ref{eq einstein y0})
follows from (\ref{eq E y1 y2 b}) and (\ref{eq einstein y i}).
Since we assume $n_1 \leq n_2$,  the two Einstein solutions of
the system satisfy  (\cite[p.313]{Wa92})
\[
y_0 \in \left [ \frac{3}{(n_1+2)^2}, 1 \right ), \quad 1 < y_2 \leq y_1 \leq 2.
\]
Hence by solving the quadratic equation (\ref{eq einstein y i}) for $y_i$ we get
\begin{equation}\label{eq y i solu}
 y_i =y_i(y_0) =  \frac{2(n_i+2)^2 y_0}{3} \left (1 -\sqrt{1- \frac{3}{y_0} \cdot \frac{1}{
(n_i+2)^2}} ~ \right),
 \end{equation}
in which we took the  negative sign in the quadratic formula because of $y_i \leq 2$.
 Hence we can define a function $\phi$ of $y_0 \in \left [ \frac{3}{(n_1+2)^2}, 1 \right )$
\begin{equation}
\phi(y_0)  \doteqdot  3(y_0-1) +2n_1 (y_1(y_0) -1) + 2n_2(y_2(y_0) -1). \label{eq phi def of y0}
\end{equation}
Note that $\phi$ is actually a convex function (see \cite[p.313]{Wa92}) which accounts for
the existence of the two Einstein points $\xi$ and $\eta$.

Let $(y_{01}, y_{11}, y_{21})$ (corresponding to $\xi$) and  $(y_{02}, y_{12}, y_{22})$
(corresponding to $\eta$) be the two solutions of (\ref{eq einstein y0}) and (\ref{eq einstein y i})
with $y_{01} < y_{02}$. Then
\begin{align}
\xi_1 & =\frac{y_{11}}{4 \Lambda_1 y_{01}} =   \frac{(n_1+2)^2 }{6  \Lambda_1} \left (1 -\sqrt{1-
\frac{3}{y_{01}} \cdot \frac{1}{(n_1+2)^2}} ~ \right)  \notag \\
& > \frac{(n_1+2)^2 }{6\Lambda_1 } \left (1 -\sqrt{1-
\frac{3}{y_{02}} \cdot \frac{1}{(n_1+2)^2}} ~ \right) =\eta_1. \label{eq xi formula}
\end{align}
Similarly we have $\eta_2 < \xi_2$. This proves $($iii$)$ and hence the lemma.
\end{proof}

\vskip .1cm
Next  we give some estimates of  $\eta$ which will be used later to study its linear stability
as a solution of (\ref{eq m=2 b 1,2 vers u}).

  \begin{lemma} \label{lem est of eta}
  The Einstein point $\eta=(\eta_1,\eta_2)$ satisfies the following properties:

 \vskip .1cm
 \noindent $($i$)$ When $n_2 \geq  n_1 \geq 2$  and  $(n_1,n_2) \neq (2,2), (2,3)$, we have
 $q_i \eta_i <   \frac{0.4661}{n_i+2}$ for $i=1,2$.

 \vskip .1cm
 \noindent $($ii$)$ When $n_1 =1$ and $n_2 \geq 2$ we have $q_1 \eta_1 < 0.1608$ and
$q_2 \eta_2 <   \frac{ 0.4661}{n_2+2}$.

 \vskip .1cm
 \noindent $($iii$)$   When $(n_1,n_2) =(2,2)$ we have $q_1 \eta_1< 0.0912 $ and
 $q_2 \eta_2 <  0.0912 $.

 \vskip .1cm
 \noindent $($iv$)$  When $(n_1,n_2) =(2,3)$ we have $ q_1 \eta_1<0.1204$ and $q_2 \eta_2
 <0.0928$.

 \vskip .1cm
 \noindent $($v$)$  When $(n_1,n_2) =(1,1)$ we have $ q_1 \eta_1<0.1303$ and $q_2 \eta_2
 <0.1303$.
\end{lemma}

\begin{proof}
We will divide the proof into five cases. In each case we find an appropriate $y_0^*$ such that
$\phi(y_0^*) <0$. Since $y_{02}$ is the larger of the two solutions of $\phi(y) =0$,
it follows that  $y_{02} > y_0^*$. Hence $y_{i2}(y_{02}) < y_{i2}(y_0^*)$ for $i=1, 2$ by (\ref{eq y i solu})
and the desired estimates will follow.

\vskip .1cm
If we choose $y_0^* = \frac{6 + \sqrt{2}}{12} $, then for all $n_i \geq 1, i=1, 2$,
\begin{equation}
 \frac{3}{y_0^*} \cdot \frac{1}{(n_i+2)^2}  < 2(\sqrt{2}-1). \label{eq y0* requirement}
\end{equation}
Note that
\begin{equation}
1-\frac{1}{2} \alpha - \frac{1}{4}\alpha ^2 < \sqrt{1-\alpha} \quad \text{for } \alpha \in (0,
2(\sqrt{2}-1)).  \label{eq alpha ineq}
\end{equation}
Applying this inequality to (\ref{eq y i solu}) with $\alpha =
\frac{3}{y_0^*} \cdot \frac{1}{(n_i+2)^2}$
 we can estimate $\phi(y_0^*)$ as follows:
\begin{align}
\phi(y_0^*) < & \,\, 3(y_0^* -1)  + \frac{3}{y_0^*} \left ( \frac{n_1}{(n_1+2)^2}
  +  \frac{n_2}{(n_2+2)^2} \right )   \label{eq est phi y0*} \\
= & \,\, 3(6- \sqrt{2}) \cdot \left (- \frac{1}{12}+ \frac{6}{17} \cdot \left ( \frac{n_1}{ (n_1+2)^2} +
\frac{n_2}{(n_2+2)^2}  \right ) \right ) \notag  \\
\leq & \,0  \,\,\notag
\end{align}
for $n_2 \geq n_1 \geq 2$ and $n_2 \geq 4$, or $(n_1,n_2) =(1,n_2)$ for any $n_2 \geq 2$.
By (\ref{eq y i solu}) and  a direct computation we also have $\phi(y_0^*) < 0$ for $(n_1,n_2) =(3,3)$.
Hence for these $(n_1,n_2)$ the corresponding
\begin{equation}
y_{02} > y_0^* =  \frac{6+  \sqrt{2}}{12}. \label{eq y02 and y*}
\end{equation}

\vskip .1cm
There are now five cases to consider.

\vskip .1cm
$($i$)$  $n_2 \geq n_1 \geq 2$ and $(n_1,n_2) \neq (2,2), (2,3)$:
By the definition of $\eta_i$, (\ref{eq y i solu}),  (\ref{eq y02 and y*}), and (\ref{eq alpha ineq}) we get
for $y_0^* = \frac{6+ \sqrt{2}}{12}$
\begin{align}
q_{i} \eta_i = & \,\,\frac{n_i+2}{6} \left (1 -\sqrt{1- \frac{3}{y_{02}} \cdot \frac{1}{
(n_i+2)^2}} ~ \right)  \notag   \\
< &  \,\, \frac{n_i+2}{6} \left (1 -\sqrt{1- \frac{3}{y^*_{0}} \cdot \frac{1}{
(n_i+2)^2}} ~ \right) \label{eq q eta i est upper}  \\
< & \,\, \frac{1}{n_i+2} \cdot \left ( \frac{1}{4 y_0^*} + \frac{3}{8(y_0^*)^2 (n_i+2)^2} \right )
\notag   \\
<& \,\, \frac{  0.4661 }{n_i+2}.  \notag
\end{align}
Note that (\ref{eq q eta i est upper}) holds whenever $y_{02} > y_0^*$.

\vskip .1cm
$($ii$)$    $n_1= 1$ and $n_2 \geq 2$:  By the definition of $\eta_i$,  (\ref{eq y02 and y*}) and
(\ref{eq q eta i est upper}) we have for $y_0^* = \frac{6 + \sqrt{2}}{12}$
\begin{align*}
& q_1 \eta_1 < \frac{1}{2} \left (1- \sqrt{1-\frac{1}{3y_0^*}} \right )< 0.1608 , \,\,\mbox{\rm and}  \\
& q_2 \eta_2 <  \frac{  0.4661 }{n_i+2}.
\end{align*}

\vskip .1cm
$($iii$)$ $(n_1,n_2) = (2,2)$: Taking $y_0^* =\frac{3+ \sqrt{2}}{6}$, we use
(\ref{eq y i solu}) to compute that
\[
 y_1(y_0^*) =y_2(y_0^*) < 1.0732.
\]
Hence we get $\phi(y_0^*) <0$ and  $y_{02} > y_{0}^* = \frac{3+ \sqrt{2}}{6}$.
We then apply (\ref{eq q eta i est upper}) to estimate
\begin{align*}
q_i \eta_i < \frac{2}{3} \left( 1- \sqrt{1 - \frac{1}{y_0^*} \cdot \frac{3}{16}} \right ) <0.0912,
\quad  i =1, 2.
\end{align*}

\vskip .1cm
$($iv$)$ $(n_1,n_2) = (2,3)$:  We take $y_0^* =\frac{10+ \sqrt{2}}{20}$. Then for this $y_0^*$ both
(\ref{eq y0* requirement}) and (\ref{eq est phi y0*})  hold for $(n_1,n_2) = (2,3)$,
hence $y_{02} > y_{0}^* = \frac{10+ \sqrt{2}}{20}$.
 By  (\ref{eq q eta i est upper}) we have
 \begin{align*}
& q_1 \eta_1 < \frac{2}{3} \left ( 1- \sqrt{1- \frac{1}{y_0^*} \cdot \frac{3}{16} } \right ) < 0.1204,   \\
& q_2 \eta_2 <  \frac{5}{6} \left ( 1- \sqrt{1- \frac{1}{y_0^*} \cdot \frac{3}{25} } \right ) < 0.0928.
 \end{align*}

\vskip .1cm
$($v$)$ $(n_1,n_2) = (1,1)$:  We take $y_0^* =\frac{3+ \sqrt{2}}{6}$. Then $y_1(y_0^*) =
y_2(y_0^*)<1.1498$ and $\phi(y_0^*)<-0.1936<0$. Hence  $y_{02} > y_{0}^* =
\frac{3+ \sqrt{2}}{6}$ and
\[
q_i \eta_i < \frac{1}{2} \left (1- \sqrt{1-\frac{1}{3y_0^*}} \right )< 0.1303.
\]
This completes the proof of the lemma.
\end{proof}

\subsection{Linearization of $\operatorname{ODE}$s (\ref{eq m=2 b 1,2 vers u})}
\label{subsec 2.4 linearlization}

We can now consider the linear stability of the vector field $(-Y_1 F_1(Y), -Y_2 F_2(Y))$ at each of
its zeros given in Lemma \ref{lem station sol QK}.

\begin{proposition} \label{prop fixed point asyp type}
For the vector field $(-Y_1 F_1(Y), -Y_2 F_2(Y))$

\vskip .1cm
\noindent $($i$)$ $(0,0)$ is a source;

\vskip .1cm
\noindent $($ii$)$ $v_1$ is  a hyperbolic point with stable eigen-direction $(1,0)$ and unstable
eigen-direction $(0,1)$;

\vskip .1cm
\noindent $($iii$)$ $v_2$ is  a hyperbolic point with unstable eigen-direction $(1,0)$ and stable
eigen-direction $(0,1)$;

\vskip .1cm
\noindent $($iv$)$ $\tilde{v}_1$ is  a source;

\vskip .1cm
\noindent $($v$)$  $\tilde{v}_2$ is  a source;

\vskip .1cm
\noindent $($vi$)$ $\eta$ is  a sink; and

\vskip .1cm
\noindent $($vii$)$  $\xi$ is a hyperbolic point.

\end{proposition}

\begin{proof}
The Jacobian of the vector field $(-Y_1 F_1(Y), -Y_2 F_2(Y))$ is given by
\begin{equation}
\mathcal{L}_Y \doteqdot  \left( \begin{array}{cc}
h_1(Y) & 8n_2q_2^2Y_1 Y_2  \\
 8n_1q_1^2Y_1 Y_2 &  h_2(Y)
\end{array} \right) ,
\end{equation}
where
\begin{align*}
& h_1(Y) \doteqdot -4(n_1 +2) q_1Y_1+3(4n_1+6)q_1^2Y_1^2 +4n_2q_2^2 Y_2^2 +\frac{1}{2}, \\
& h_2(Y) \doteqdot -4(n_2 +2) q_2Y_2+3(4n_2+6)q_2^2Y_2^2 +4n_1q_1^2 Y_1^2 +\frac{1}{2}.
\end{align*}

\vskip .1cm
$($i$)$ At $Y=(0,0)$ we have
\begin{equation}
\mathcal{L}_{(0,0)}  =  \left( \begin{array}{cc}
\frac{1}{2} & 0  \\
 0 & \frac{1}{2}
\end{array} \right).  \label{eq L at 0,0 linear op}
\end{equation}
Hence $(0,0)$ is a source.

\vskip .1cm
 $($ii$)$  At $Y=v_1$ we have
\[
\mathcal{L}_{v_1}  =  \left( \begin{array}{cc}
- \frac{n_1+1}{2n_1+3} & 0  \\
 0 & \frac{n_1}{(2n_1+3)^2} + \frac{1}{2}
\end{array} \right).
\]
Hence $v_1$ is a hyperbolic point. It is clear that $(1,0)$ is the stable eigen-direction and
$(0,1)$ is the unstable eigen-direction.

\vskip .1cm
 $($iii$)$ This follows from its symmetry to  case $($ii$)$.

\vskip .1cm
 $($iv$)$  At $Y=\tilde{v}_1$ we have
\[
\mathcal{L}_{\tilde{v}_1}  =  \left( \begin{array}{cc}
n_1+1 & 0  \\
 0 & n_1+ \frac{1}{2}
\end{array} \right).
\]
Hence $\tilde{v}_1$ is a  source.

\vskip .1cm
 $($v$)$ This follows from its symmetry to  case $($iv$)$.

\vskip .1cm
 $($vi$)$ At  $Y=\eta$,  using $F_1(\eta) =F_2(\eta) =0$, we have
\[
h_1(\eta) = (8n_1+6)q_1^2\eta_1^2 - E(\eta) ,  \quad
h_2(\eta) = (8n_2+6)q_2^2\eta_2^2 - E(\eta),
\]
and $\mathcal{L}_{\eta} = - E(\eta) I  + \beta_\eta $ where for $\chi =(\chi_1,\chi_2)$
\begin{equation}
\beta_\chi \doteqdot \left( \begin{array}{cc}
(8n_1+6)q_1^2\chi_1^2  & 8n_2q_2^2 \chi_1 \chi_2  \\
  8n_1q_1^2 \chi_1 \chi_2 &  (8n_2+6)q_2^2\chi_2^2
\end{array} \right).   \label{eq beta chi def}
\end{equation}
We need to analyze the eigenvalues of $\beta_\eta$.

\vskip .1cm
Suppose $\chi_1>0$ and $\chi_2>0$. Using the diagonal matrix $ [\chi_1, \chi_2]$,
we compute that
\begin{equation}
\alpha_\chi \doteqdot [\chi_1^{-1}, \chi_2^{-1}] \cdot \beta_\chi \cdot [\chi_1, \chi_2]
= \left( \begin{array}{cc}
(8n_1+6)q_1^2\chi_1^2  & 8n_2q_2^2  \chi_2^2  \\
  8n_1q_1^2 \chi_1^2 &  (8n_2+6)q_2^2\chi_2^2
\end{array} \right).  \label{eq alpha eta def}
\end{equation}
Note that since $\det \alpha_\chi >0$, we may assume that $\rho_1(\chi) \geq \rho_2(\chi) >0$ are
the two eigenvalues of $\alpha_\chi$. Then the eigenvalues of $\mathcal{L}_{\chi}$
are given by $-E(\chi) + \rho_i(\chi)$.

A simple calculation of the eigenvalues
of $\alpha_\eta$ using the Einstein condition gives
\begin{equation}
-E(\eta)+ \rho_1(\eta) = 3q_1^2 \eta_1^2 +  3q_2^2 \eta_2^2 -\frac{1}{2}
+\sqrt{A(\eta) } ,  \label{eq eta matric eigen}
\end{equation}
where
\begin{align*}
A(\chi) \doteqdot &  (4n_1+3)^2q_1^4\chi_1^4 + (4n_2+3)^2q_2^4\chi_2^4 \\
& + (32n_1 n_2-24n_1 -24n_2 -18) q_1^2 q_2^2 \chi_1^2 \chi_2^2 .
\end{align*}
It is easy to see that
\[
A(\eta) <\left ( (4n_1+3)q_1^2\eta_1^2 + (4n_2+3)q_2^2\eta_2^2 \right  )^2.
\]
So for all $n_2 \geq n_1 \geq 1$
\begin{equation}
 -E(\eta)+ \rho_1(\eta) <  ( 4n_1+6)q_1^2\eta_1^2 + (4n_2+6)q_2^2\eta_2^2 -\frac{1}{2} .
 \label{eq eigen value eta 1}
\end{equation}

Hence by Lemma \ref{lem est of eta}(i)  we have
\begin{equation*}
  -E(\eta)+ \rho_1(\eta)
<  (4n_1+6)  \left ( \frac{0.4661}{n_1+2} \right  )^2
+ (4n_2+6) \left ( \frac{0.4661 }{n_2+2} \right )^2 -\frac{1}{2}    < 0
\end{equation*}
when  $n_2 \geq n_1 \geq 2$ and $(n_1,n_2) \neq (2,2), (2,3)$.
By Lemma \ref{lem est of eta}(ii) we have for $n_1=1$ and  $n_2 \geq 2$
\[
 -E(\eta)+ \rho_1(\eta)   < 10 \cdot  \left ( 0.1608 \right  )^2 + (4n_2+6) \cdot
 \left ( \frac{0.4661 }{n_2+2} \right )^2 -\frac{1}{2} <0.
 \]
When $(n_1,n_2) =(2,2)$,  Lemma \ref{lem est of eta}(iii)
we have
 \[
  -E(\eta)+ \rho_1(\eta)   < 14 \cdot \left ( 0.0912\right )^2 +  14 \cdot
  \left ( 0.0912\right )^2 -\frac{1}{2} <0.
  \]
When  $(n_1,n_2) =(2,3)$, by Lemma \ref{lem est of eta}(iv)  we have
 \[
  -E(\eta)+ \rho_1(\eta)   < 14 \cdot \left ( 0.1204 \right )^2
   +18 \cdot \left ( 0.0928 \right )^2 -\frac{1}{2} <0.
  \]
  When  $(n_1,n_2) =(1,1)$, by Lemma \ref{lem est of eta}(v)  we have
 \[
  -E(\eta)+ \rho_1(\eta) < 10 \cdot (0.1303)^2+ 10 \cdot (0.1303)^2 -\frac{1}{2} <0.
  \]
 We have therefore proved that  $-E(\eta)+ \rho_1(\eta)   < 0$ for all $(n_1,n_2)$.
 Hence $\eta$ is a sink.

\vskip .1cm
 $($vii$)$   Similar to case $($vi$)$ we have  $\mathcal{L}_{\xi} = - E(\xi) I  + \beta_\xi $
 where $\beta_\xi$ is obtained from (\ref{eq beta chi def}) by setting $\chi=\xi$.
For the corresponding matrix $\alpha_\xi$ defined by (\ref{eq alpha eta def}) with eigenvalues
 $\rho_1(\xi) \geq \rho_2(\xi) >0$, we obtain  $-E(\xi) +\rho_1(\xi) \geq -E(\xi) + \rho_2(\xi)$
for the two eigenvalues of $\mathcal{L}_\xi$.

To see that $ -E(\xi)+ \rho_2(\xi)  <0$, it is easy to check that for all $n_1, n_2, q_1\xi_1, q_2 \xi_2$,
\[
A(\xi) > (3q_1^2 \xi_1^2 +  3q_2^2 \xi_2^2)^2,
\]
from which it follows that
\[
-E(\xi)+ \rho_2(\xi) = 3q_1^2 \xi_1^2 +  3q_2^2 \xi_2^2 -\frac{1}{2}
-\sqrt{A(\xi)}  <0.
\]
Note that in the line above we took the negative square root in front of
$\sqrt{A(\xi)}$ since $\rho_2(\xi)$ is by choice the smaller eigenvalue.

We will give an indirect proof that $-E(\xi)+ \rho_1(\xi) >0$ when we establish
 the second claim of Proposition \ref{Prop eternal solution behavior xi hyperbol}(ii).
  Hence $\xi$ is a hyperbolic point.
In view of this the proof of Proposition \ref{prop fixed point asyp type} is complete.
\end{proof}

\section{\bf Ancient solutions on $\R\PP^3$ fibre bundles} \label{sec 3 ancient sol existence metric behav}

In this section we first prove the existence of ancient solutions on the
$\R\PP^3$ fibre bundles $\bar{P}$ over a product of two QK manifolds. Then we turn to
study the asymptotic behaviors of the metric tensors.

\subsection{ Existence of ancient solutions on $\R\PP^3$ fibre bundles} \label{subsec existence on RP3 bundle}

We shall consider two regions in $\mathbb{R}^2$: a convex region
$$\Omega_1 \doteqdot \{Y: \, F_1(Y) \geq 0 \text{ and } F_2(Y) \geq 0 \},$$
and a (non-convex) region
$$\Omega_2 \doteqdot  \,\, \mbox{\rm  bounded component of }\{Y: \, Y_1 \geq 0, \, Y_2 \geq 0,
 \,F_1(Y) \leq 0, \text{ and } F_2(Y) \leq 0 \}.$$

\begin{theorem} \label{thm psi b one two}
Suppose $m=2$ and the QK manifolds $(M_i^{4n_i}, g_i)$ have positive scalar curvature,
 i.e., $\Lambda_i >0$ for $i=1, 2$. Then the $\operatorname{ODE}$ system
$($\ref{eq m=2 a evol m=2 vers 2}$)$ and $($\ref{eq m=2 b 1,2}$)$
has the following types of long-time solutions $(\psi(\tau),b_1(\tau),b_2(\tau))$
with $0 \leq \tau < + \infty$ corresponding to the solutions $Y(u)$ of
 $($\ref{eq m=2 b 1,2 vers u}$)$ :

\vskip .1cm
\noindent $($i$)$  There is a continuous $1$-parameter family of  solutions $Y(u)$
 with $\lim_{u \rightarrow \infty} Y(u) =\eta$ and a corresponding
 continuous $1$-parameter family of ancient solutions  $\bar{g}(\tau)$
  of the Ricci flow on $\bar{P}$. Furthermore, as special cases,

\vskip .1cm
\noindent $($ia$)$ 	if the initial data $Y(0) =(\frac{\psi(0)}{b_1(0)}, \frac{\psi(0)}{b_2(0)})  \in \Omega_1$,
then the solution $Y(u)$ stays in $\Omega_1$ and  $\lim_{u \rightarrow \infty} Y(u) =\eta$;

\vskip .1cm
 \noindent $($ib$)$ if the initial data $Y(0)   =(\frac{\psi(0)}{b_1(0)},
 \frac{\psi(0)}{b_2(0)})   \in \Omega_2$, then the solution  $Y(u)$ stays in $\Omega_2$
  $\lim_{u \rightarrow \infty} Y(u) =\eta$.

\vskip .1cm
 \noindent $($ii$)$ There are two distinct solutions $Y(u)$ with
 $\lim_{u \rightarrow \infty} Y(u) =\xi$ and correspondingly two distinct
 ancient solutions $\bar{g}(\tau)$  on $\bar{P}$.
\end{theorem}

\begin{proof}
$($i$)$ By Proposition \ref{prop fixed point asyp type}(vi) the local stable manifold at $\eta$
has dimension $2$, so the existence of a continuous $1$-parameter family of ancient solutions $Y(u)$
 follows from the Hartman-Grobman Theorem.

From (\ref{eq m=2 a evol m=2 vers u})
  we can conclude that $\psi(u)$ exists for all $u \in [0,\infty)$
and that $\psi(u)^{-1}$ is bounded from above.
Hence from  (\ref{eq def u via psi}) we get $\tau (u) \rightarrow \infty$ when $u \rightarrow \infty$.
In view of the discussion after (\ref{eq m=2 b 1,2 vers u}), we get a continuous
 $1$-parameter family of ancient solutions $(\psi(\tau),b_1(\tau),b_2(\tau))$.
This in turn gives a  continuous $1$-parameter family of ancient solutions  $\bar{g}(\tau)$
(modulo time translation and parabolic scaling).

\vskip .1cm

$($ia$)$ Consider a solution $Y(u)$ with $Y(0)
 \in \Omega_1 \setminus \{ \xi,\eta \}$.
If at some time $u_0 \geq 0$ we have $Y(u_0) \in \partial  \Omega_1 \setminus \{ \xi,\eta \}$,
then either  $F_1(Y(u_0)) =0$ or $F_2(Y(u_0)) =0$. In the first case,
computing at $Y(u_0)$, we obtain
$$\nabla F_1 = (2 (n_1 +2)q_1 - 2 (4n_1+6)  q_1^2 Y_1, -8n_2 q_2^2 Y_2), $$
so that
$$\nabla F_1 \cdot (0, -Y_2 F_2 (Y)) = 8n_2 q_2^2 Y_2^2 F_2(Y) > 0,$$
where the last inequality follows from the assumption that $Y(u_0) \neq \xi$ or $\eta$.
Similarly, if  $F_2(Y(u_0)) =0$, we have $\nabla F_2 \cdot (-Y_1 F_1(Y),0) > 0$.
Hence the solution $Y(u)$ stays in $\Omega_1$ for all $u$.
The property  $\lim_{u \rightarrow \infty} Y(u) =\eta$ will be proved in Theorem
\ref{Prop eternal solution behavior xi hyperbol} below.

\vskip .1cm

$($ib$)$  The vector field  $ -(Y_1 F_1(Y),Y_2 F_2(Y))$  has 4 zeros $\{0,\eta,v_1,v_2\}$
in $\Omega_2$. Take a solution $Y(u)$ of  (\ref{eq m=2 b 1,2 vers u}) with
$Y(0) \in \Omega_2 \setminus \{ 0,\eta,v_1,v_2 \}$.
If there is a time $u_0$ at which the solution $Y(u)$ hits the boundary where $F_1(Y) =0$,
then computing at $Y(u_0)$ as above one gets
\[
\nabla F_1 \cdot (0, -Y_2 F_2 (Y)) = 8n_2 q_2^2 Y_2^2 F_2(Y) \leq 0,
\]
with equality iff we are at the fixed points $\eta$ or $v_1$.
If the solution $Y(u)$ hits the boundary where $F_2(Y) =0$  at time $u_0$,
we have at $Y(u_0)$
\[
\nabla F_2 \cdot (-Y_1 F_1 (Y), 0) = 8n_1 q_1^2 Y_1^2 F_1(Y) \leq 0,
\]
with equality iff we are at the fixed points $\eta$ or $v_2$.

If the solution $Y(u)$ hits the part of the boundary where $Y_2 =0$ at time
$u_0$, then  we have $F_1(Y(u_0)) \leq 0$
and  $F_2(Y(u_0)) < 0$, so that
\[
(0,1)  \cdot (-Y_1F_1(Y), -Y_2 F_2(Y)) = 0.
\]
Finally, if the solution $Y(u)$ hits the boundary where $Y_1 =0$ at time $u_0$,
then  we have $F_1(Y(u_0)) < 0$ and  $F_2(Y(u_0)) \leq 0$. Hence
\[
(1, 0)  \cdot (-Y_1F_1(Y), -Y_2 F_2(Y)) = 0.
\]
Combining the above analysis we conclude that $Y(u)$ remains in $\Omega_2$
for all $u \geq 0$. The property  $\lim_{u \rightarrow \infty} Y(u) =\eta$ will be
proved in Theorem \ref{Prop eternal solution behavior Omega 2} below.

\vskip .1cm
$($ii$)$ The existence of two distinct solutions $Y(u)$  follows from the Hartman-Grobman Theorem
since the local stable manifold at $\xi$ has dimension $1$.
  By analogous arguments to those for case (i) we obtain
  the existence of two ancient solutions $\bar{g}(\tau)$.
This proves case (ii) and hence also completes the proof of the Theorem.
\end{proof}

\begin{rmk} \label{pseudo-Rflow}
We may also consider the cases where $\Lambda_1 <0, \, Y_1<0$ and/or
$\Lambda_2 <0, \, Y_2<0$. These cases correspond to pseudo-Riemannian
Ricci flows. Notice that in (\ref{eq m=2 b 1,2}), if $q_i < 0$ and $b_i < 0$,
we can multiply the equation by $-1$ and obtain the equation for $|q_i|$
and $|b_i|$. Therefore, we can deduce the existence of ancient solutions
from that of the positive case.

Another way to get  a pseudo-Riemannian flow  is to assume that the signs of
$\tau$ and $\psi$ in equations (\ref{eq m=2 a evol m=2 vers 2})
and (\ref{eq m=2 b 1,2}) are both negative.
By changing the sign of $\tau$ and $\psi$ and replacing $q_i$ by $|q_i|$
we  obtain the equation for $|\tau|$, $|\psi|$ and $|q_i|$.
Again we can deduce from the existence of ancient solutions for the positive case that there are
 immortal solutions of pseudo-Riemannian metrics on $\bar{P}$ which are
negative definite on the fibres and positive definite on the base.

Note that there are many non-symmetric homogeneous negative QK manifolds,
see e.g. \cite{Co96}. As well, C. Lebrun showed that the moduli space of complete
QK structures on $\R^{4n}$ is infinite-dimensional \cite{Le91}. Finally, we mention the
recent thesis of M. Dyckmanns \cite{Dy15} in which he constructed complete negative QK manifolds
 that are not locally homogeneous for all dimensions $4n$.
\end{rmk}

\subsection{Limiting behavior of the ancient solutions} \label{subsec 3.3 limiting behavior}

Next we consider the longtime behavior of the solutions in Theorem \ref{thm psi b one two}.
The arguments are analogous to those for the torus bundles in \cite[\S 3.5]{LW16},
 but because the fibres now have positive curvature,
 some of the conclusions are different.

\begin{theorem} \label{Prop eternal solution behavior xi hyperbol}
Let $Y(u)$ be one of the ancient solutions in Theorem \ref{thm psi b one two}$($ia$)$ with initial value
$Y(0) \in \Omega_1 \setminus \{\xi \}$. Then

\noindent $($i$)$ the forward limit $\lim_{u \rightarrow \infty} Y(u) =\eta$, and

\noindent $($ii$)$ there is exactly one solution whose backwards limit  $\lim_{u \rightarrow -\infty} Y(u)
= \xi$.
\end{theorem}

\begin{proof}
We use the notations in the proof of Proposition \ref{prop fixed point asyp type}.

\vskip .1cm
 $($i$)$ Let $\omega_Y$ be the the  $\omega$-limit set of the flow line $Y(u)$.
Since
\[
\frac{d}{du} E(Y(u)) = - 8n_1q_1^2 Y_1^2 F_1(Y) - 8n_2q_2^2 Y_2^2 F_2(Y) \leq 0,
\]
with equality if and only if $F_1(Y) =F_2(Y) =0$, it follows that $\omega_Y =\{\eta \}$ or $\{ \xi \}$.
Below we rule out the possibility that $\omega_Y =\{\xi \}$.

We \textbf{claim} that for any $Y(0) \in \Omega_1 \setminus \{ \xi \}$ there is an $i_0$
such that $Y_{i_0}(0) < \xi_{i_0}$. Assuming this claim, we compute that $\frac{d Y_{i_0}}{du}(u)
=-Y_{i_0}(u) F_{i_0}(u) \leq 0$ and so $\lim_{u \rightarrow \infty} Y_{i_0}(u)< \xi_{i_0}$.
Hence  $\omega_Y =\{\eta \}$ must hold.

We next prove the above claim by considering two cases. We rewrite the equations $F_1(Y) =F_2(Y) =0$
in the $(q_1Y_1,q_2Y_2)$-plane as loci of ellipses:
\begin{subequations}
\begin{align}
& \frac{ \left (q_1Y_1 - \frac{n_1+2}{4n_1+6} \right )^2 }{ \left ( \frac{n_1 +1}{4n_1+6} \right )^2 }
+\frac{(q_2Y_2)^2}{  \frac{(n_1+1)^2}{ (4n_1 +6) \cdot 4n_2} } =1,  \label{eq ellipse 1} \\
&\frac{(q_1Y_1)^2}{  \frac{(n_2+1)^2}{ (4n_2 +6) \cdot 4n_1} }  +
\frac{ \left (q_2Y_2 - \frac{n_2+2}{4n_2+6} \right )^2 }{ \left ( \frac{n_2 +1}{4n_2+6} \right )^2 }
 =1.  \label{eq ellipse 2}
\end{align}
\end{subequations}

$($ia$)$ $n_2 \geq n_1 \geq 1$ and $(n_1,n_2) \neq (1,1)$. From (\ref{eq ellipse 1}) we conclude that
\begin{equation}
q_2 \xi_2 \leq \sqrt{ \frac{(n_1+1)^2}{ (4n_1 +6) \cdot 4n_2}}.  \label{eq tem q2xi 2 neg}
\end{equation}
From (\ref{eq ellipse 2}) we get
\begin{equation}
q_2 \xi_2 =\frac{n_2 +2 \pm \sqrt{(n_2+1)^2 -(4 n_2 +6) \cdot 4 n_1 q_1^2 \xi_1^2} }{4n_2 +6}.
\label{eq q 2 xi 2 expression}
\end{equation}
Since for  $n_2 \geq n_1 \geq 1$ and $(n_1,n_2) \neq (1,1)$ one can check that
\[
\frac{n_2 +2}{4n_2+6} >  \sqrt{ \frac{(n_1+1)^2}{ (4n_1 +6) \cdot 4n_2}},
\]
it follows from (\ref{eq tem q2xi 2 neg}) that we need to take negative sign in (\ref{eq q 2 xi 2 expression}),
and the square root in the equation cannot be zero. We have proven that  $(q_1 \xi_1, q_2 \xi_2)$
lies in the lower-right quarter of the ellipse  in the $(q_1Y_1,q_2Y_2)$-plane defined by $F_2(Y) =0$  and
cannot be the right vertex of the ellipse.
This implies that every point $(Y_1,Y_2) \in \Omega_1\setminus \{ \xi \}$ satisfies $Y_1 < \xi_1$.

$($ib$)$ $n_2 = n_1 = 1$. Then by  (\ref{eq ellipse 1}) and  (\ref{eq ellipse 2})
we know that the set
\[
\{ (q_1Y_1,q_2Y_2): \,\,(Y_1,Y_2) \in \Omega_1 \}
\]
is a convex set
in $(q_1Y_1,q_2Y_2)$-plane and is symmetric with respect to the diagonal.
This implies that every point $(Y_1,Y_2) \in \Omega_1\setminus \{ \xi \}$ satisfies either $Y_1 < \xi_1$
or $Y_2 < \xi_2$, otherwise $\xi$ will be an interior point of $\Omega_1$, which is not true.
Now the claim is proved.

\vskip .1cm
 $($ii$)$  Let
 \begin{align*}
 & a_i \doteqdot 2 (n_i+2)q_i - (4n_i +6) \cdot 2 q_i^2 \xi_i , \\
 & b_i \doteqdot 8n_i q_i^2 \xi_i, \quad i=1, 2.
 \end{align*}
 We \textbf{claim} the following estimate
 \begin{equation}
 \frac{b_1}{a_1+b_1} + \frac{b_2}{a_2+b_2} >1.  \label{eq a i bi est key}
 \end{equation}

To verify the claim, note that $\nabla F_1(\xi) =(a_1, -b_2)$ and $\nabla F_2(\xi) =( -b_1, a_2)$.
From (ia) and (ib) above, there is a vector $(-\delta_1, -\delta_2)$ with $\delta_i >0, i=1, 2$,
lying in the interior of the tangent cone $T_{\xi} \Omega_1$, i.e.,
there is a solution $(\lambda_1, \lambda_2)$ with $\lambda_i >0, i=1, 2$, to the following
equations,
\begin{subequations}
\begin{align}
 a_1 \lambda_1 -b_1 \lambda_2 & = - \delta_1, \label{eq a i bi lam 1}  \\
-b_2 \lambda_1 +  a_2 \lambda_2 & = - \delta_2,  \label{eq a i bi lam 2}  \\
\lambda_1 + \lambda_2 & = 1 .  \label{eq a i bi lam 3}
\end{align}
\end{subequations}

Note that it follows from $q_i \xi_i < \frac{1}{2}$ that
\[
a_i +b_i = 2q_i \left ( (n_i +2) -6q_i \xi_i \right ) > 0, \quad i=1, 2.
\]
We compute using (\ref{eq a i bi lam 1}), (\ref{eq a i bi lam 2}) and (\ref{eq a i bi lam 3}) that
\begin{align*}
\frac{b_1}{a_1+b_1} + \frac{b_2}{a_2+b_2} &=\frac{\lambda_1(a_1+b_1)+ \delta_1}{a_1+b_1} +
\frac{\lambda_2(a_2+b_2)+ \delta_2}{a_2+b_2} \\
& = 1+ \frac{ \delta_1}{a_1+b_1} + \frac{ \delta_2}{a_2+b_2} ,
\end{align*}
which gives the claim (\ref{eq a i bi est key}).

Next we \textbf{claim} that the eigenvalue $-E(\xi) +\rho_1(\xi)$ of $\mathcal{L}_{\xi}$, defined in
Proposition \ref{prop fixed point asyp type}(vii), is positive. To see this, we compute
\begin{align*}
0< & \,\,\frac{b_1}{a_1+b_1} + \frac{b_2}{a_2+b_2} -1 \\
 =& \,\,\frac{8n_1q_1^2 \xi_1^2}{E(\xi) -6 q_1^2\xi_1^2}
+ \frac{8n_2q_2^2 \xi_2^2}{E(\xi) -6 q_2^2\xi_2^2} -1      \\
= & \,\,\frac{G(\xi, E(\xi))}{(E(\xi) -6 q_1^2\xi_1^2) \cdot (E(\xi) -6 q_2^2\xi_2^2)  },
\end{align*}
where for $ \tilde{E} \in \mathbb{R}$
\begin{align*}
G(\xi, \tilde{E}) \doteqdot & -\tilde{E}^2 + \left ( (8n_1+6)q_1^2 \xi_1^2+ (8n_2+6)q_2^2 \xi_2^2
 \right ) \tilde{E}  \\
&  - (48n_1 +48n_2 +36)q_1^2 q_2^2 \xi_1^2 \xi_2^2 .
\end{align*}
It is easy to check that $E(\xi) -6 q_i^2\xi_i^2 >0, \, i=1, 2$. Hence we conclude that
$G(\xi, E(\xi)) >0$.

 Since $\rho_1(\xi)$ is an eigenvalue of the matrix $\alpha_{\xi}$ defined by (\ref{eq alpha eta def}),
 $\rho_1(\xi)$ satisfies $G(\xi, \rho_1(\xi)) =0$.
Hence
\begin{align}
0<& \, G(\xi, E(\xi)) -G(\xi, \rho_1(\xi)) \notag  \\
=& (-E(\xi) +\rho_1(\xi)) \cdot \left ( E(\xi) + \rho_1(\xi) - (8n_1+6)q_1^2 \xi_1^2
-(8n_2 +6) q_2^2 \xi_2^2 \right ).  \label{eq E rho 1 difference tem}
\end{align}
Note that $E(\xi) -6 q_1^2 \xi_1^2 >0$, and that by applying the Perron-Frobenius theory
to the matrix $\alpha_{\xi}$ we have $ \rho_1(\xi) \geq 8n_1q_1^2 \xi_1^2
+ (8n_2 +6) q_2^2 \xi_2^2 $ (see  p.76 of \cite{Ga59} for example). Hence the second factor in
(\ref{eq E rho 1 difference tem}) is positive and so  $-E(\xi) +\rho_1(\xi) >0$.

Now we can finish the proof of (ii) by using the Hartman-Grobman Theorem
since the local unstable manifold at $\xi$ has dimension $1$.
\end{proof}

\vskip .1cm
When the solution $Y(u)$ starts in $\Omega_2$, its longtime behavior is given by

\begin{theorem} \label{Prop eternal solution behavior Omega 2}
Let $Y(u)$ be one of the ancient solutions in Theorem \ref{thm psi b one two}$($ib$)$ with initial value
$Y(0) \in \Omega_2 \setminus (\overline{0 v_1} \cup \overline{0 v_2} \cup \{ \eta \} )$
where $\overline{0 v_i}$ denotes the line segment joining $0$ and $v_i$. Then

\noindent $($i$)$ the forward limit $\lim_{u \rightarrow \infty} Y(u) =\eta$.

\noindent $($ii$)$  There is exactly one solution whose backwards limit  $\lim_{u \rightarrow -\infty} Y(u)
= v_1$. We denote the corresponding flow line by $\gamma_1$.

\noindent $($iii$)$ There is exactly one solution whose backwards limit  $\lim_{u \rightarrow -\infty}
Y(u)  = v_2$. We denote the corresponding flow line by $\gamma_2$.

\noindent $($iv$)$ Let $\Omega_2^* \subset \Omega_2$ be the closed region bounded by
$\overline{0 v_1}$,  $\overline{0 v_2}$,  $\gamma_1$ and $\gamma_2 \cup \{\eta \}$.
For each $Y(0)$ in the interior of $\Omega_2^*$ the solution $Y(u)$ has  backwards limit
$\lim_{u \rightarrow -\infty} Y(u) = 0$.
\end{theorem}

\begin{proof}
We use the notations in the proof of Proposition \ref{prop fixed point asyp type}.

\vskip .1cm
 $($i$)$ Let $\omega_Y$ be the $\omega$-limit set of the flow line $Y(u)$.
Then
\[
\frac{d}{du} E(Y(u)) = - 8n_1q_1^2 Y_1^2 F_1(Y) - 8n_2q_2^2 Y_2^2 F_2(Y) \geq 0
\]
with equality if and only if $Y_1F_1(Y) =Y_2 F_2(Y)=0$. As $\omega_Y$ is flow-invariant
and connected, it follows that $\omega_Y =\{ 0 \}, \, \{v_1\}, \,
\{ v_2 \}$, or $\{\eta \}$.
Below we rule out the possibilities that $\omega_Y =\{ 0 \}, \, \{v_1\}$, or  $\{ v_2 \}$.

From (ia) and(ib) in the proof of Theorem \ref{Prop eternal solution behavior xi hyperbol}
we conclude that $(q_1 \eta_1, q_2 \eta_2)$ lies in the right-lower quarter of the ellipse  in
$(q_1Y_1,q_2Y_2)$-plane defined by $F_2(Y) =0$  and cannot be the right vertex of the ellipse.
Likewise, the point lies in the upper-left quarter of the ellipse  $F_1(Y) =0$.
These facts imply that that $Y_{i}(0) < \eta_{i}$ for any $Y(0) \in  \Omega_2 \setminus \{\eta \}$.
But  $\frac{d Y_{i}}{du}(u) =-Y_{i}(u) F_{i}(u) \geq 0$, and the inequality is  strict
for at least one $i$ by the assumption on $Y(0)$. This implies that  $\omega_Y =\{\eta \}$.

\vskip .1cm
$($ii$)$ and $($iii$)$ Both follow from the Hartman-Grobman Theorem as in the proof
of Theorem \ref{Prop eternal solution behavior xi hyperbol}(ii).

\vskip .1cm
$($iv$)$ Since the boundary of  $\Omega_2^*$ consist of flow lines, by the uniqueness of solutions
of ODEs, $Y(u)$ must stay in $\Omega_2^*$.  From $\frac{d Y_{i}}{du}(u)>0$ proved in (i), the backward limit
must be one of $0$, $v_1$, and $v_2$. The last two limits are ruled out by the uniqueness in (ii) and (iii).
\end{proof}

\subsection{Asymptotic behavior of metric tensors of the ancient solutions}

 We begin with the following.
\begin{theorem}\label{thm growth ai bi QK Omega 1}
Let  $(\psi(\tau), b_1(\tau), b_2(\tau))$ be one of the ancient solutions of the $\operatorname{ODE}$
system $($\ref{eq m=2 a evol m=2 vers 2}$)$ and $($\ref{eq m=2 b 1,2}$)$ given in
Theorem \ref{thm psi b one two}$($i$)$.
Then we have the following estimates and asymptotics.

\vskip .1cm
\noindent $($i$)$ The domain for $\tau$ contains $[0, \infty)$. For any $\varepsilon >0$ small
there is a $\tau_{\varepsilon} >0$ such that for any $\tau > \tau_{\varepsilon}$ and for $i=1,2$,
\begin{align}
 & (E(\eta) -\varepsilon) ( \tau -\tau_{\varepsilon} )  \leq \psi(\tau) - \psi(\tau_{\varepsilon} )
\leq \left( E(\eta) + \varepsilon\right) (\tau -\tau_{\varepsilon} ) ,
 \label{eq est psi QK} \\
& -\varepsilon  (\tau -\tau_{\varepsilon} )  \leq b_i(\tau) -
 b_i(\tau_{\varepsilon} ) - \left ( 2(n_i +2)
q_i - 6q_i^2 \eta_i \right )    (\tau - \tau_{\varepsilon} )
  \leq \varepsilon (\tau - \tau_{\varepsilon} ).
\label{eq est bi QK}
\end{align}
Furthermore, we have
\[
\lim_{\tau \rightarrow \infty} \frac{\psi (\tau)}{b_i(\tau)} = \eta_i, \quad
\lim_{\tau \rightarrow \infty} \frac{ \psi (\tau)}{\tau} = E(\eta),
\]
and so the rescaled metrics $\frac{1}{\tau} \bar{g}(\tau)$ converge
in the Gromov-Hausdorff topology to a multiple of the Einstein metric corresponding to $\eta$ as
$\tau$ approaches $\infty$.

\vskip .1cm
\noindent $($ii$)$ For the solution in Theorem \ref{Prop eternal solution behavior xi hyperbol}$($ii$)$
there is a finite $T_1>0$ such that when $ u\rightarrow -\infty$ \,the corresponding
 $\tau \rightarrow -T_1^+$. We have
\[
\lim_{\tau \rightarrow - T_1^+} \psi(\tau) =0, \quad \lim_{\tau \rightarrow - T_1^+} b_i(\tau) =0,
\quad \lim_{\tau \rightarrow - T_1^+} \frac{\psi (\tau)}{b_i(\tau)} = \xi_i, \quad
\lim_{\tau \rightarrow -T_1^+} \frac{ \psi (\tau)}{T_1 +\tau} = E(\xi).
\]
It follows  that geometrically, as $\tau \rightarrow - T_{1}^+$, the
$\R\PP^3$ fibre bundle  $\bar{P}$, equipped with the metric $\bar{g} (\tau)$,
collapses to a point. Furthermore, the rescaled metrics
$\frac{1}{ T_1+\tau } \, \bar{g} (\tau)$ converges to a multiple of the
Einstein metric corresponding to $\xi$  as $\tau \rightarrow -T_1^+$.
\end{theorem}

\begin{proof}
\noindent $($i$)$  We first note that the range of $\tau(u)$ (as determined by
(\ref{eq def u via psi})) contains $[0, \infty)$ for any of the solutions $Y(u)$ under consideration.
The inequalities in (\ref{eq est psi QK}) follow from $\frac{d \psi}{d\tau} = E(Y)$ and $\lim_{u
\rightarrow \infty} Y(u) =\eta$ and those in (\ref{eq est bi QK}) follow from
(\ref{eq m=2 b 1,2}) and $\lim_{u \rightarrow \infty}  Y(u) =\eta$.
Note that  $2(n_i +2) q_i - 6q_i^2 \eta_i >0$, since we have $q_i \eta_i < \frac{1}{2}$
from the proof  of Lemma \ref{lem station sol QK}(ii).
The values of the limits follow from the convergence of $Y(u)$ and
 (\ref{eq m=2 a evol m=2 vers 2}).

\vskip .1cm
 $($ii$)$ The proof is similar to that of  \cite[Theorem 3.1.3(ii)]{LW16}.
We first show that $\lim_{u \rightarrow - \infty} \tau(u) = -T_1 $ for some finite
positive $T_1$. Since $\lim_{u \rightarrow -\infty} \, Y(u) = \xi$, we have
$$ \frac{d \psi}{d \tau} > \frac{1}{2} + n_1 q_1^2 \xi_1^2 + n_2 q_2^2 \xi_2^2$$
for $\tau$ corresponding to $u \in (-\infty, u_*]$. Here $u_*$ is a sufficiently negative number.
Integrating this inequality over the $\tau$-interval corresponding to $[u, u_*]$ gives
$$0 < \tau(u_*) - \tau(u)  < \psi(u_*) \left(\frac{1}{2} + n_1 q_1^2 \xi_1^2 + n_2 q_2^2
\xi_2^2 \right)^{-1}. $$
This implies the existence of finite $T_1 >0$.

Next one shows that $\lim_{\tau \rightarrow -T_1^+} \psi(\tau) =0$. The limit certainly
exists as $\psi(\tau)$ is positive and increasing in $\tau$. If the limit is positive,
we immediately get a contradiction to $\lim_{\tau \rightarrow -T_1^+} \, u(\tau) = \infty$
upon integrating $\frac{du}{d\tau} = \frac{1}{\psi}$.
We omit the remaining details as they are straight-forward.
\end{proof}

\vskip .1cm
When the solution $Y(u)$ starts from $\Omega_2$ the asymptotic behavior of the metric tensor
is given by the next theorem.

\begin{theorem}\label{thm growth ai bi QK Omega 2}
Let $(\psi(\tau), b_1(\tau), b_2(\tau))$ be a solution  of the $\operatorname{ODE}$ system
$($\ref{eq m=2 a evol m=2 vers 2}$)$ and $($\ref{eq m=2 b 1,2}$)$.

\vskip .1cm
\noindent $($i$)$  If  it corresponds to the solution in Theorem
\ref{Prop eternal solution behavior Omega 2}$($ii$)$ then
there is a finite $T_2>0$ such that when $ u\rightarrow -\infty$ the corresponding $\tau \rightarrow
-T_2^+$.  We in addition have
\begin{align*}
& \lim_{\tau \rightarrow - T_2^+} \psi(\tau) =0, \quad \lim_{\tau \rightarrow - T_2^+} b_1(\tau) =0,
\quad \lim_{\tau \rightarrow - T_2^+} b_2(\tau) > 0,  \\
& \lim_{\tau \rightarrow - T_2^+} \frac{\psi (\tau)}{b_1(\tau)} =\frac{1}{(4n_1+6)q_1}, \quad
\lim_{\tau \rightarrow - T_2^+} \frac{\psi (\tau)}{b_2(\tau)} =0, \quad
\lim_{\tau \rightarrow -T_2^+} \frac{\psi (\tau)}{T_2 +\tau} = E(v_1).
\end{align*}
Geometrically,  as $\tau \rightarrow - T_2^+$, the $\R\PP^3$ fibre bundle  $\bar{P}$, equipped with the
metric $\bar{g} (\tau)$, collapses along the fibres and the first factor
to some QK metric on the second factor $M_2$. The rescaled metrics $\frac{1}{ T_2+ \tau }\,
\bar{g} (\tau)$
converge, as $\tau \rightarrow -T_2^+$, to product $P_1 \times \R^{4n_2}$ in the Gromov-Hausdorff topology,
where $P_1$ is equipped with a multiple of the ``squashed" Einstein metric $($see the remark below$)$
 and $\R^{4n_2}$ is equipped with the Euclidean metric.

\vskip .1cm
\noindent $($ii$)$  If it corresponds to the solution in
Theorem \ref{Prop eternal solution behavior Omega 2}$($iii$)$ then
there is a finite $T_3>0$ such that when $ u\rightarrow -\infty$ the corresponding $\tau \rightarrow
-T_3^+$.  We further have
\begin{align*}
& \lim_{\tau \rightarrow - T_3^+} \psi(\tau) =0, \quad \lim_{\tau \rightarrow - T_3^+} b_1(\tau) > 0,
\quad  \lim_{\tau \rightarrow - T_3^+} b_2 (\tau) =0,  \\
& \lim_{\tau \rightarrow - T_3^+} \frac{\psi (\tau)}{b_1(\tau)} =0, \quad
 \lim_{\tau \rightarrow - T_3^+} \frac{\psi (\tau)}{b_2(\tau)} =\frac{1}{(4n_2 +6)q_2}, \quad
\lim_{\tau \rightarrow -T_3^+} \frac{\psi (\tau)}{T_3 +\tau} = E(v_2).
\end{align*}
Geometrically,  as $\tau \rightarrow - T_3^+$, the bundle  $\bar{P}$, equipped with the
metric $\bar{g} (\tau)$, collapses along the fibres and the second factor
to some QK metric on the first factor $M_1$. The rescaled metrics $\frac{1}{ T_3 + \tau }\,
\bar{g} (\tau)$
converge, as $\tau \rightarrow -T_3^+$, to product $\R^{4n_1}\times P_2$ in the Gromov-Hausdorff topology,
where $\R^{4n_1}$ is equipped with the Euclidean metric
and  $P_2$ is equipped with a multiple of the ``squashed" Einstein metric.

\vskip .1cm
\noindent $($iii$)$  If it corresponds to a solution in Theorem
\ref{Prop eternal solution behavior Omega 2}$($iv$)$ then there is a finite
$T_4>0$ such that when $ u\rightarrow -\infty$ the corresponding $\tau \rightarrow
-T_4^+$.  We also have
\begin{align*}
 \lim_{\tau \rightarrow - T_4^+} \psi(\tau) =0, \quad \lim_{\tau \rightarrow - T_4^+} b_i(\tau) > 0,
, \quad \lim_{\tau \rightarrow - T_4^+} \frac{\psi (\tau)}{b_i(\tau)} =0, \quad
\lim_{\tau \rightarrow -T_4^+} \frac{\psi (\tau)}{T_4 +\tau} = \frac{1}{2}.
\end{align*}
Geometrically,  as $\tau \rightarrow - T_4^+$, the bundle  $\bar{P}$, equipped with the
metric $\bar{g} (\tau)$, collapses along the fibres
to the product metric
\[
 \lim_{\tau \rightarrow - T_4^+} b_1(\tau) \cdot g_1 +  \lim_{\tau \rightarrow - T_4^+} b_2(\tau)
  \cdot g_2
\]
on $M_1 \times M_2$. By choosing different initial data the ratio $\frac{ \lim_{\tau \rightarrow -
T_4^+} b_2(\tau)}{ \lim_{\tau \rightarrow - T_4^+} b_1(\tau)}$ can be any positive number.
The rescaled metrics $\frac{1}{T_4 + \tau }\,\bar{g} (\tau)$
converge, as $\tau \rightarrow -T_4^+$, to product $\R\PP^3 \times \R^{4(n_1 +n_2)}$ in
the Gromov-Hausdorff topology,
where $\R\PP^3$ is equipped with the metric with constant sectional  curvature $1/4$
and $\R^{4(n_1 +n_2)}$ is equipped with the Euclidean metric.
\end{theorem}

\begin{proof}
$($i$)$  The proof is similar to that of Theorem \ref{thm growth ai bi QK Omega 1}$($ii$)$
except for the assertion that $\lim_{\tau \rightarrow - T_2^+} b_2(\tau) > 0$.
To prove this, we first give an estimate of $T_2$. Using again (\ref{eq m=2 a evol m=2 vers u})
and the monotonicity of $E(Y(u))$, we  have $E(0) \leq E(Y(u)) =\frac{d \ln \psi}{du} \leq E(\eta)$.
Integrating this inequality over the interval $[u,0]$ we get
\[
\psi(0)\, e^{E(\eta)u} \leq \psi(\tau(u)) =\frac{d \tau}{du} \leq \psi(0)\, e^{E(0)u},
\]
which in particular gives $\lim_{u \rightarrow -\infty} \, \psi(\tau(u)) = 0$.
Integrating again, we get for $u<0$
\[
\frac{\psi(0)}{E(\eta)} \left (1 -e^{E(\eta)u} \right ) \leq - \tau(u) \leq
\frac{\psi(0)}{E(0)} \left (1 -e^{E(0)u} \right ).
\]
Hence by letting $u \rightarrow -\infty$ we get
\begin{equation}
\frac{\psi(0)}{E(\eta)} \leq T_2 \leq \frac{\psi(0)}{E(0)}  =2 \psi(0).  \label{eq est T3 time}
\end{equation}

By integrating (\ref{eq m=2 b 1,2}) we get  that
$b_2(\tau) \geq b_2(0) + 2q_2 (n_2 + 2) \tau $ for  $\tau<0$.
By using (\ref{eq est T3 time}) we have
\[
\lim_{\tau \rightarrow -T_2^+} b_2(\tau) \geq b_2(0) - 4q_2 (n_2 + 2) \psi(0)
= b_2(0)(1 -  4q_2 (n_2 + 2) Y_2(0) ) >0.
\]

To get the last inequality,  we may need to shift the initial time so that $q_2Y_2(0) < (4(n_2+2))^{-1}$,
which is always possible.
 The remaining limits follow from the above and L'Hopital's rule.

\vskip .1cm
 $($ii$)$  Since the proof is similar to that of $($i$)$ above, we omit it.

\vskip .1cm
 $($iii$)$ Except for the last limit  $\frac{ \lim_{\tau \rightarrow -
T_4^+} b_2(\tau)}{ \lim_{\tau \rightarrow - T_4^+} b_1(\tau)}$,
 the rest of the proof is similar to that
 of $($i$)$ above and so we omit it. To deduce the limit, we use (\ref{eq L at 0,0 linear op})
 together with Sternberg's Theorem (\cite[Theorem 4]{St57}) on the $C^k$-conjugacy of solutions
 near a fixed point of a nonlinear $\operatorname{ODE}$ system satisfying non-resonance conditions to
 solutions of its linearization.  It follows that there is a diffeomorphism $W$ from a neighborhood of the origin
 in $(Y_1, Y_2)$-space to some neighborhood of $0 \in \R^2$ such that the flow lines $Y(u)$ satisfy
 $$ Y(u) = (W^{-1} \circ (k_1e^{u/2}, k_2 e^{u/2}) \circ \, W)(Y(0)) $$
 where $u$ is sufficiently close to $-\infty$. Taking the $u$-derivative of this equation and letting $u$
 tend to $-\infty$ we see that the tangent map $W^{-1}_{*}$ maps directions at $0$ with slope
 $k_2/k_1$ to directions at $0$ in $Y$-plane. In particular, as $\overline{0v_1}$ and $\overline{0v_2}$
 are themselves flow lines, it follows that there are flow lines $Y(u)$ emanating from $0$
 which lie in between $\overline{0 v_1}$ and  $\overline{0 v_2}$ with any positive limiting
slope.
\end{proof}

\begin{rmk} \label{squashed-metrics}
Here we clarify the term ``squashed" Einstein metric for readers not familiar with Einstein metrics
on the principal ${\rm SO}(3)$ bundle $P$ over a compact QK manifold. We will use the
notation in \S \ref{subsec 2.1 intro to QK}. It is well-known (see \cite[Proposition 14.85]{Be87})
that $P$ admits two Einstein metrics of connection type. Indeed, they can be found by solving
the quadratic equation resulting from (\ref{eq E y1 y2 b}) by dropping the terms with index $2$.
The two solutions are then given by $Y_1 = \frac{1}{2q_1}$ and $Y_1 = \frac{1}{(4n_1 + 6)q_1}$.
The latter is what we refer to as the ``squashed" Einstein metric. The terminology comes from the
case when the QK manifold is quaternionic projective space $\HH\PP^n$ and $P$
is $S^{4n+3}$. The ``squashed" Einstein metric is the Jensen metric, obtained from the constant
curvature metric by uniformly shrinking the Hopf fibres. It is of interest to physicists since
the dimension of the space of Killing spinors is $1$. There are also pseudo-Riemannian analogues
of these Einstein metrics (see \cite[Remark 14.86b]{Be87}), which appear as limits in our construction of
pseudo-Riemannian flows.
\end{rmk}

\begin{theorem}\label{thm growth ai bi QK Omega 3}
Let $(\psi(\tau), b_1(\tau), b_2(\tau))$ be one of the ancient solutions of the $\operatorname{ODE}$
system $($\ref{eq m=2 a evol m=2 vers 2}$)$ and $($\ref{eq m=2 b 1,2}$)$ in Theorem
\ref{thm psi b one two}$($ii$)$.
Then the domain for $\tau$ contains $[0, \infty)$ and for any small $\varepsilon >0$ there is a
 $\tau_{\varepsilon}>0$ such that for $\tau \geq \tau_{\varepsilon}$ and for $i=1, 2$,
\begin{align}
&  (E(\xi) -\varepsilon) (\tau -\tau_{\varepsilon} ) \leq \psi(\tau) - \psi(\tau_{\varepsilon}) \leq
(E(\xi) + \varepsilon) (\tau -\tau_{\varepsilon} ) ,     \label{eq est psi QK iii} \\
& (2n_i+1)\,q_i ( \tau  -\tau_{\varepsilon}) \leq b_i(\tau) - b_i(\tau_{\varepsilon}) \leq 2(n_i +2)\, q_i
( \tau -\tau_{\varepsilon}).
\label{eq est bi QK iii}
\end{align}
Furthermore, we have
\[
\lim_{\tau \rightarrow \infty} \frac{\psi (\tau)}{b_i(\tau)} = \xi_i, \quad
\lim_{\tau \rightarrow \infty} \frac{ \psi (\tau)}{\tau} = E(\xi),
\]
and so the rescaled metrics $\frac{1}{\tau} \bar{g} (\tau)$ converge
in the Gromov-Hausdorff topology to a multiple of the Einstein metric corresponding to $\xi$ as
$\tau$ approaches $\infty$.
\end{theorem}

\begin{proof}
Inequalities (\ref{eq est psi QK iii}) follows from $\frac{d \psi}{d\tau} = E(Y)$ and $\lim_{u
\rightarrow \infty} Y(u) =\xi$.
Since $q_i \xi_i < \frac{1}{2}$, we may choose $\tau_{\varepsilon}$ so that
by using $($\ref{eq m=2 b 1,2}$)$ we have  $(2n_i +1) q_i \leq \frac{d b_i}{d \tau} \leq 2(n_i +2) q_i $
for $\tau \geq \tau_{\varepsilon}$.
Then the inequality (\ref{eq est bi QK iii}) follows.
The remaining convergences  follow easily.
\end{proof}

With the asymptotics given in Theorems \ref{thm growth ai bi QK Omega 1},
\ref{thm growth ai bi QK Omega 2} and \ref{thm growth ai bi QK Omega 3},
one easily deduces the following
\begin{corollary} \label{cor noncollapsing}
 $($i$)$ Each of the ancient solutions $\bar{g} (\tau), \, \tau \in [0,\infty)$,
in Theorem \ref{thm psi b one two},
is $\kappa$-noncollapsed at all scales for some $\kappa > 0$.

\noindent{$($ii$)$ } Each of the solutions $\bar{g}(\tau)$ on $ \tau \in (-T_i, 0],
1 \leq i \leq 4$ in  Theorems \ref{thm growth ai bi QK Omega 1}$($ii$)$ and
\ref{thm growth ai bi QK Omega 2}$($i$)$-$($iii$)$, is $\kappa$-noncollapsed at all scales
for some $\kappa > 0$.
\end{corollary}

\section{\bf $C^0$ Cheeger-Gromov convergence implies $C^\infty$ Cheeger-Gromov convergence
and Type I solutions of the Ricci flow } \label{subsec 4 QK curva prop}

Although we can compute the curvature tensor of the ancient solutions we found in
Theorem \ref{thm psi b one two} as we did for the torus bundle case
in \cite{LW16} and conclude that these solutions are of type I,
in this section we will instead prove the solutions are of type I by an indirect method.
This relies on a result relating $C^0$ and $C^{\infty}$ convergence of Ricci flow solutions
under suitable hypotheses. We believe that this result may of be independent interest
and have other applications.

Recall the following definition of  convergence in pointed $C^k$ Cheeger-Gromov  topology,
$0 \leq k \leq \infty$.

\begin{definition} \label{def name 4.1}
Let $I$ be a fixed interval. A sequence $\left\{  \left( {M}_{k}^{n}, g_{k}( t ), p_{k}\right)
 \right\}, \,t \in I$, of complete pointed solutions of the Ricci flow converges in the $C^k$
 (pointed) Cheeger-Gromov  topology
  to a complete pointed solution $\left( {M}_{\infty}^{n}, g_{\infty} (  t),
  p_{\infty}\right), \,  t \in I$, of the Ricci flow if there exist

\vskip .1cm
 \noindent $($i$)$ an exhaustion $\left\{  U_{k}\right\}  _{k\in\mathbb{N}}$ of
${M}_{\infty}$ by open sets with $p_{\infty}\in U_{k}$ for each $k$, and

\vskip .1cm
 \noindent $($ii$)$  a sequence of diffeomorphisms $\Phi_{k}:U_{k}\rightarrow
 \Phi_{k}\left(  U_{k}\right)  \subset {M}_{k}$ with $\Phi
_{k}\left(  p_{\infty}\right)  =p_{k}$

\noindent such that $\left(  U_{k},\Phi_{k}^{\ast}\left ( \left.  g_{k}\left(
t\right)  \right\vert _{\Phi_k(U_{k})}\right )  \right) $ converges in the $C^{k}$ topology to
$\left( {M}_{\infty},g_{\infty}\left(  t\right)  \right)  $ uniformly
on any compact subset in ${M}_{\infty}\times I$.
\end{definition}

We have the following strengthening of $C^0$ Cheeger-Gromov convergence of the Ricci flow.
\begin{theorem}\label{thm 1 no choice}
Let  $\{ (M^n_k,  g_k(t), p_k) \}, \, t \in [0,1]$, be a sequence of complete pointed
 solutions of the Ricci flow with bounded curvature for each $k$.
We assume that

\vskip .1cm
\noindent $($a$)$ the scalar curvature $R_{g_k}(x,t) \geq -R_*$ on $M_k \times[0,1]$ for
 some constant $R_* > 0$,

\vskip .1cm
\noindent $($b $)$ the sequence $\{ (M_k,  g_k(t), p_k) \}, \, t \in [0,1]$,
converges in the $C^0$ Cheeger-Gromov topology
 to complete pointed solution $(M^n_\infty,  g_{\infty}(t), p_{\infty}), \, t \in [0,1]$, of the Ricci flow.

\vskip .1cm
\noindent Then the sequence $\{ (M_k,  g_k(t), p_k) \}, \, t \in (0,1]$, subconverges
 to $(M_\infty,  g_{\infty}(t), p_{\infty}), \, t \in (0,1]$, in $C^\infty$
 Cheeger-Gromov topology.
\end{theorem}

\begin{proof}
Fix $\alpha =1$ and let $\delta $ and $\epsilon_0 <1$ be the two positive constants in Perelman's
pseudolocality theorem  (\cite[Theorem 10.1]{Pe02}).
We \textbf{claim} that  for any compact subset $K \subset M_{\infty}$
  there is a $r_* \in (0, \min\{1/2, R_*^{-1/2} \})$ such that for any $x_* \in K$ and any 
  $t \in [0,1]$ the ball  $B_{g_\infty(t)}(x_*, 2r_*)$ is $\delta/2$-almost
isoperimetrically Euclidean in the sense
\begin{equation}
(\operatorname{Area}_{ g_{\infty}(t)}(\partial \Omega))^n \geq \left(1 - \frac{\delta}{2} \right) c_n
( \operatorname{Vol}_{g_\infty(t)}(\Omega))^{n-1} \label{eq delta e almost euc isoper}
\end{equation}
for any regular domain $\Omega \subset B_{g_\infty(t)}(x_*, 2 r_*)$,
where $c_n =n^n \omega_n$ is the Euclidean isoperimetric constant.

{\noindent \it Proof of the claim.} Because every Riemannian metric is infinitesimally a 
Euclidean metric and the isoperimetric property is stable for  nearly equivalent metrics,
given any $(\tilde{x}, \tilde{t}) \in M_{\infty} \times [0,1]$ there are  $\tilde{r}>0$ and $\tilde{\eta} >0$
sufficiently small such that
\begin{equation}
(\operatorname{Area}_{ g_{\infty}(t)}(\partial \Omega))^n \geq \left(1 - \frac{ \delta}{2} \right) c_n
( \operatorname{Vol}_{g_\infty(t)}(\Omega))^{n-1}  \label{eq 4r almost euc isoper}
\end{equation}
for any regular domain $\Omega \subset B_{g_\infty(t)}(\tilde{x}, 4 \tilde{r})$ and any $t \in
I_{\tilde{t},\tilde{\eta}}$. We define interval
$$  I_{\tilde{t},\tilde{\eta}} \doteqdot
\begin{cases}
(\max\{ \tilde{t} - \tilde{\eta}, 0 \}, \min \{\tilde{t} + \tilde{\eta},1 \}) & \mbox{if }
\tilde{t} - \tilde{\eta}> 0 \text{ and } \tilde{t} + \tilde{\eta} <1 \\
[0, \min \{\tilde{t} + \tilde{\eta},1 \}) & \mbox{if } \tilde{t} - \tilde{\eta} \leq 0  \text{ and } \tilde{t}
+ \tilde{\eta} <1 \\
(\max\{ \tilde{t} - \tilde{\eta}, 0 \},1] & \mbox{if }  \tilde{t} - \tilde{\eta}> 0  \text{ and } \tilde{t}
 + \tilde{\eta} \geq 1 \\
[0,1] & \mbox{if } \tilde{t} - \tilde{\eta} \leq 0  \text{ and } \tilde{t} + \tilde{\eta} \geq 1
\end{cases}.
$$
Note that
$$ \{ \cup_{t \in I_{\tilde{t},\tilde{\eta}}} B_{g_\infty(t)}(\tilde{x}, 2 \tilde{r}) \times \{t \}  ,
\, (\tilde{x}, \tilde{t}) \in M_{\infty} \times [0,1] \}
$$
is an open cover of $M_{\infty} \times [0,1]$.
Hence  for any compact subset $K \subset M_{\infty}$, we can extract a finite subcovering
$ \{ \cup_{t \in I_{\tilde{t}_i,\tilde{\eta}_i}} B_{g_\infty(t)}(\tilde{x}_i, 2 \tilde{r}_i) \times \{ t \}
, \, i =1, 2, \cdots, i_K \} $ of $K \times [0,1]$.
We choose $r_*$ to be any positive number less than both $\min_{i=1, 2, \cdots, i_K } \tilde{r}_i$ and $1/2$.

Now we verify the conclusion of the claim. Given any $(x_*,t_*) \in K \times [0,1]$ there is an $i$ such that
$(x_*,t_*) \in \cup_{t \in I_{\tilde{t}_i,\tilde{\eta}_i}} B_{g_\infty(t)}(\tilde{x}_i, 2 \tilde{r}_i) \times \{ t \}$.
Then given any regular domain $\Omega \subset B_{g_\infty(t_*)}(x_*, 2 r_*)$ we have
$ \Omega  \subset  B_{g_\infty(t_*)}(\tilde{x}_i, 4 \tilde{r}_i)$ where $t_* \in  I_{\tilde{t}_i,\tilde{\eta}_i}$.
Hence by the choice of $\tilde{x}_i$, $\tilde{t}_i$, $\tilde{r}_i$,
 and $\tilde{\eta}_i$ in (\ref{eq 4r almost euc isoper}) we know that
 inequality (\ref{eq delta e almost euc isoper}) holds
for  $\Omega $ and the claim holds.

\smallskip

Since limit manifold $M_{\infty}$ may be noncompact, below we will use balls centered at
$p_{\infty}$ with larger and larger radius to provide the exhaustion in Definition \ref{def name 4.1}.
Fix any $D>0$. By assumption (b) the metrics $\Phi_k^*   g_{k} ( t )$  converges in $C^0$ topology to
$g_{\infty}(t)$ uniformly  on $\cup_{t \in [0,1]} B_{g_{\infty}(t)} (p_{\infty}, D+2) \times \{ t\}$.
Using this fact, the claim, and the property that the isoperimetric property is stable for
nearly equivalent metrics, we conclude that there is an $r_* \in (0, \min\{1/2, R_*^{-1/2} \})$ 
which is independent of $k$  and  $k_0$ such that for any $k \geq k_0$, $t\in [0,1]$,
and $x_0 \in B_{g_k(t)}(p_k,D)$ we have
\begin{equation}
(\operatorname{Area}_{ g_{k}(t)}(\partial \Omega))^n \geq (1 - \delta) c_n
( \operatorname{Vol}_{g_k(t)}(\Omega))^{n-1} \label{eq almost euc isoper seq k}
\end{equation}
for any regular domain $\Omega \subset B_{g_k(t)}(x_0, r_*)$.

By (\ref{eq almost euc isoper seq k}), $r_* \leq R_*^{-1/2}$, and by the assumption (a)
we may apply Perelman's
pseudolocality theorem to the Ricci flow $(M_k, \left . g_k(t) \right |_{t \in [t_*,t_* +
(\epsilon_0 r_*)^2]})$ with $k\geq k_0$ and
$t_* \in [0,1 - (\epsilon_0 r_*)^2]$.
We then obtain  the interior  curvature estimate
\begin{equation}
|\operatorname{Rm}_{g_k}|(x, t) \leq \frac{1}{t-t_*}+ \frac{1}{(\epsilon_0 r_*)^2}
\label{eq perel est seq}
\end{equation}
where $d_{g_k(t)}(x, x_0) < \epsilon_0 r_*$, $x_0 \in B_{g_k(t)}(p_k,D)$, and
$t  \in (t_*, t_* + (\epsilon_0 r_*)^2]$.

Given any $D >0$ and $\bar{t}>0$, by applying inequality (\ref{eq perel est seq}) for some $t_* \in
[\max \{0,\bar{t} - (\epsilon_0 r_*)^2 \}, 1- (\epsilon_0 r_*)^2 ]$
we get the uniform curvature bound
\begin{equation}
|\operatorname{Rm}_{g_k}|(x, t) \leq \frac{1}{\bar{t}} + \frac{2}{(\epsilon_0 r_*)^2}
\label{eq bdd radius bdd curv}
\end{equation}
for $k \geq k_0$, $x \in B_{g_k(t)}(p_k,D)$,  and $t  \in [\bar{t},1]$.

\vskip .1cm
Next we are going to prove that the sequence of injectivity radii $\operatorname{inj}_{g_k(1/2)}(p_k)$
have an uniform lower bound independent of $k$. By the assumption (b) we know that
\[
\lim_{k \rightarrow \infty} \operatorname{Vol}_{g_k(1/2)}(B_{g_k(1/2)}(p_k, 1)) =
 \operatorname{Vol}_{g_{\infty}(1/2)}(B_{g_{\infty}(1/2)}(p_\infty, 1)).
 \]
By (\ref{eq bdd radius bdd curv}) there is a constant $K_0$ such that
the curvature $|\operatorname{Rm}_{g_k}|(x, 1/2) \leq K_0$ for any $x \in B_{g_k(1/2)}(p_k, 1)$
and $k \geq k_0$. Hence by Lemma \ref{thm CGT local version sep 20, 2016} below
  there is a constant $\iota_0>0$ independent of $k$
 such that $\operatorname{inj}_{g_k(1/2)}(p_k) \geq \iota_0$.

Now we state a local version of a theorem of Cheeger, Gromov, and Taylor  (see \cite[Theorem 4.7]{CGT82}). 

\begin{lemma}\label{thm CGT local version sep 20, 2016}
 For any positive constants $r_{0},  k_0, v_0$,  and dimension $n$, there exists a constant $\iota_{0}>0$
 depending only on $r_0, k_0,v_0, n$
such that if for some precompact ball $B(p, r_0)$  in a Riemannian manifold $(M^{n},g)$
the sectional curvature $\vert \operatorname{sect} (x) \vert \leq k_0$ for all $x \in B(p,r_0)$  and the volume
$\operatorname{Vol} (B(p,r_0)) \geq v_0$, then the injectivity radius
$\operatorname{inj} (p)  \geq \iota_{0}$.
\end{lemma}

\begin{proof}
It is well known that $\operatorname{inj}(p) = \min \{ \text{conjugate radius at } p,
\frac{1}{2} \ell(p) \}$, where $\ell(p)$ denotes the length of the shortest (nontrivial) closed geodesic 
starting and ending at $p$. The conjugate radius at $p$ is bounded from below by a constant $\iota_{1}
\doteqdot \min \{r_0, \frac{\pi}{\sqrt{k_0}} \}$.
It follows from \cite[Theorem 4.3]{CGT82} that $\ell(p)$ is bounded from below by a constant $\iota_{2}>0$
 depending only on $r_0, k_0, v_0, n$. The theorem now is proved.
\end{proof}

\vskip .1cm
With the curvature bound (\ref{eq bdd radius bdd curv}) and the injectivity radius lower bound,
we can now apply Hamilton's local compactness theorem (\cite[Theorem 16.1]{Ha93}) to sequence
  $\{ (M_k,  g_k(t), p_k) \}, \, t   \in (0,1]$, with base time $t =1/2$ to conclude its subconvergence
 to $(\tilde{M}^n_\infty, \tilde{g}_{\infty}(t), \tilde{p}_{\infty}), \, t \in (0,1]$, in $C^\infty$
 Cheeger-Gromov topology. By the uniqueness of the $C^0$ Cheeger-Gromov limit we
 know that $(\tilde{M}_\infty, \tilde{g}_{\infty}(t), \tilde{p}_{\infty})$ is isometric to
 $({M}_\infty, {g}_{\infty}(t), {p}_{\infty})$ for  each $t \in (0,1]$.
\end{proof}

The following corollary shows that the ancient solutions we constructed in this article
are of type I as time approaches to $-\infty$.

\begin{corollary} \label{cor used to get type I asypm}
Let $ \bar{g}(\tau) \doteqdot \bar{g}_{\vec{a}(\tau), \vec{b}(\tau)}, \, \tau \in [0, \infty)$,
be an ancient solution of the backwards Ricci flow on the compact fibre bundles $\bar{P}$ as defined in
\S \ref{subsec 5.2 metric anciet QK}.
Assume that limits  $ \lim_{\tau \rightarrow \infty} \frac{a_i(\tau)}{\tau} \doteqdot \hat{a}_i \neq 0$
 and $\lim_{\tau \rightarrow \infty} \frac{b_i(\tau)}{\tau}  \doteqdot \hat{b}_i \neq 0$
 for $i=1, \cdots, m$.
Then the solution $ \bar{g}(\tau)$ is of type I when $\tau \rightarrow \infty$.
\end{corollary}

\begin{proof}
It suffices to show that for any $\tau_k \rightarrow \infty$ there is a constant $C$ independent of $k$
such that  $ \sup_{p \in \bar{P}} \tau_k |\operatorname{Rm}_{\bar{g}(\tau_k)}(p) | \leq C$ for
some subsequence of $\tau_k$.
Note that $(\bar{P},  \bar{g}_{\vec{\hat{a}}, \vec{\hat{b}}})$ is a compact
 Riemannian manifold where $\vec{\hat{a}} =(\hat{a}_1, \cdots, \hat{a}_m)$ and
 $\vec{\hat{b}} =(\hat{b}_1, \cdots, \hat{b}_m)$.

Define the sequence of Ricci flow solutions $g_k(t) =(\tau_k)^{-1} \bar{g}(\frac{3}{2}\tau_k -\tau_k t),
\, t \in [0,1]$.
Since
\[
 \lim_{k \rightarrow \infty} \frac{a_i(\frac{3}{2}\tau_k -\tau_k t)}{\tau_k} =\left(\frac{3}{2}
  -t \right) \hat{a}_i,
 \quad \lim_{k \rightarrow \infty} \frac{b_i(\frac{3}{2}\tau_k -\tau_k t)}{\tau_k}
  =\left(\frac{3}{2} -t \right) \hat{b}_i,
 \]
the sequence of solutions  $\{ g_k(t)\}$ satisfies the assumptions of Theorem \ref{thm 1 no choice}
and it follows that  $g_k(\frac{1}{2}) = (\tau_k)^{-1} \bar{g}(\tau_k) $ subconverges in $C^\infty$
Cheeger-Gromov topology to the compact
Riemannian manifold $(\bar{P},  \bar{g}_{\vec{\hat{a}},  \vec{\hat{b}}})$.
Since $\tau_k |\operatorname{Rm}_{\bar{g}(\tau_k)} |(p) =
 |\operatorname{Rm}_{g_k(\frac{1}{2})} |(p)$ for any $p \in M$,
we have the subconvergence of 
$\sup_{p \in \bar{P}} \tau_k |\operatorname{Rm}_{\bar{g}(\tau_k)} |(p)$
to $\sup_{p \in \bar{P}} |\operatorname{Rm}_{\bar{g}_{\vec{\hat{a}},  \vec{\hat{b}}}}|(p)$.
Hence the subsequence of  $\{ \tau_k |\operatorname{Rm}_{\bar{g}(\tau_k)} |(p) \}$ is bounded.
\end{proof}

\begin{theorem} \label{thm anc sol curva prop}
The ancient solutions  $\bar{g}_{\vec{a} (\tau), \vec{b} (\tau)}$ of the Ricci flow on $\bar{P}$ with $m=2$
in  Theorem  \ref{thm psi b one two} are of Type I when $\tau \rightarrow \infty$.
The Ricci curvature of these solutions  are
positive  definite for large enough $\tau$.
The Ricci curvature of those solutions in  Theorem  \ref{thm psi b one two}$($ia$)$ and $($ib$)$  are
positive for each $\tau$.
\end{theorem}

\begin{proof}
The first assertion follows from Corollary \ref{cor used to get type I asypm} directly.

For the second assertion, the formulas (\ref{eq quaternion b}) and (\ref{eq Rc tilde U m=2})
imply that the Ricci tensors are positive definite provided that $Y_i < \frac{n_i +2}{6q_i}$ for $ i=1, 2$.
Since $q_i \eta_i < q_i \xi_i < \frac{1}{2}$, we conclude that any solution $Y(u)$ in  Theorem
\ref{thm psi b one two} satisfies the condition above
 when $u$ is large enough.

For the third assertion,  in fact in the regions $\Omega_i, i=1, 2$, we actually have $Y_i < \frac{1}{2q_i}$.
\end{proof}

\section{\bf Examples of ancient solutions  with $m\geq 3$ }
 \label{subsec 5 QK m geq 3 same structure}

In this section we shall consider a special case of the
Riemannian fibre bundles $(\bar{P}, \bar{g}_{\vec{a},\vec{b}})$ in
\S \ref{sec 2, ODE and linearizaion} with $m \geq 3$ in order to
exhibit continuous families of non-collapsed ancient solutions on even-dimensional
non-K\"ahler manifolds. For the special case we assume that the QK manifolds
$(M_i^{4n_i}, g_i)$ are positive and have the same dimension, i.e., $n_1 =\cdots =n_m \doteqdot d$.
For convenience we shall assume that the constants $\Lambda_i =n_i+2=d+2$ for all $i$,
so that the constants $q_i =1$.

 \subsection{Properties of Einstein metrics on $\bar{P}$ with $m \geq 3$}
Recall that for a solution $\bar{g}_{\vec{a}(\tau),\vec{b}(\tau)}=\bar{g}(\tau)$ of the backwards Ricci flow
equations (\ref{eq quaternion a evol}) and (\ref{eq quaternion b}) we had set
$\hat{a} \doteqdot \sum_{k=1}^m a_k$.  We now  let $X_k \doteqdot \frac{a_k}{\hat{a}}$
and $Y_k \doteqdot \frac{\hat{a}}{b_k}$, so that we have the constraint $\sum_{k=1}^m X_k =1$
for the variables $X_k$. Let $E(X,Y)  \doteqdot \frac{m}{2} +1+ \sum_{k=1}^m 4d X_k^2 Y_k^2$,
so that by (\ref{eq d a hat d tau}) we have $\frac{d \hat{a}}{ d \tau} = E(X,Y)$.
Analogous to the $m=2$ case we define a new independent variable $u$ by
$u = \int_0^{\tau} \frac{1}{\hat{a}(\zeta)} d \zeta$. Then it follows from
(\ref{eq quaternion a evol}) and (\ref{eq quaternion b}) that
\begin{subequations}
\begin{align}
&\frac{d X_k}{d u} = \frac{1}{2}+ X_k +4dX_k^2 Y_k^2 -X_k E(X,Y), \label{eq sec 5 ode a} \\
& \frac{d Y_k}{d u} = -Y_k \left ( 2(d+2)Y_k -6X_k(1-X_k)Y_k^2 -E(X,Y) \right ).
\label{eq sec 5 ode b}
\end{align}
\end{subequations}
It is easy to check that condition $\sum_{k=1}^m X_k =1$ is preserved by (\ref{eq sec 5 ode a}).
Since the metrics we shall construct must satisfy this constraint, we will only consider
solutions of the above system which satisfy $\sum_{k=1}^m X_k =1$ as well as the positivity
conditions $X_k >0$ and $Y_k >0$.
Given a solution of (\ref{eq sec 5 ode a}) and (\ref{eq sec 5 ode b}) we can find
$\hat{a}(u)$ from $\frac{d \ln \hat{a}}{du} =E(X,Y)$ first, picking up an integration
constant in the process. Then we find $\tau(u)$ from $\frac{d \tau}{du} =
\hat{a}(u)$, and finally recover $a_k(\tau)$ and $b_k(\tau)$ from $X_k$ and $Y_k$.
Hence we get solutions of the backwards Ricci flow depending on the integration constant.
Note that the net effect of the integration constant is a parabolic rescaling of the
solutions of the backwards Ricci flow. Hence given a solution of (\ref{eq sec 5 ode a})
and (\ref{eq sec 5 ode b}) we get only one solution $\bar{g}(\tau)$ of the backwards
Ricci flow (modulo time translation and parabolic rescaling).

One can easily determine all the fixed points of the $\operatorname{ODE}$ system
(\ref{eq sec 5 ode a}) and (\ref{eq sec 5 ode b}). Here we will examine those fixed
points with all $Y_k \neq 0$. This leads to the equations
\begin{subequations}
\begin{align}
&\frac{1}{2}+ X_k +4dX_k^2 Y_k^2 -X_k E(X,Y) =0,  \label{eq eins sec 5 a} \\
& 2(d+2)Y_k -6X_k(1-X_k)Y_k^2 -E(X,Y) =0.  \label{eq eins sec 5 b}
\end{align}
\end{subequations}
Let $E(X,Y) \doteqdot 2 \Lambda \hat{a}$. Then after some simple calculation we obtain
\begin{align*}
&\frac{1}{4}+ \frac{a_k}{2\hat{a}} +2d \frac{a_k^2}{ b_k^2} = \Lambda a_k , \\
& d+2 -3 \left (1- \frac{a_k}{\hat{a}} \right ) \frac{a_k}{b_k} =\Lambda b_k.
\end{align*}
By the Ricci curvature formulas in \S \ref{subsec 5.2 metric anciet QK}, these
are exactly the Einstein condition for the metric determined by $a_k, b_k$.
Hence we have proved the first part of the following
\begin{lemma}\label{lem 5.1 called}
$($i$)$ The zeros of the vector field given by $($\ref{eq sec 5 ode a}$)$ and $($\ref{eq sec 5 ode b}$)$  with
$Y_k \neq 0$ for each $k$ correspond to Einstein metrics $g^{\pm}$ on the bundle $\bar{P}$.

\vskip .1cm
\noindent $($ii$)$ Corresponding to these Einstein metrics we have $X_k =\frac{1}{m}$ and
 $Y_k = \frac{m}{\sqrt{d}}\, \beta^{\pm}$ where
 \[
\beta^{\pm} = \frac{(d+2) \left (1 \pm \sqrt{1- 2(1+2m^{-1}) d (d+2)^{-2} -3(1-m^{-1})(1+2 m^{-1})
(d+2)^{-2}  } \right ) }{2 \sqrt{d} \left (2 +3(1- m^{-1}) d^{-1} \right )} .
 \]
Note that $\beta^{\pm}$  are the solutions of quadratic equation
 \[
 \left (2+\frac{3(m-1)}{md} \right ) \beta ^2 - \frac{d+2}{\sqrt{d}}\, \beta +\frac{m+2}{4m} =0.
 \]
\end{lemma}

\begin{proof}
$($ii$)$  Since $E(X,Y) = \frac{m}{2}+1+4m (\beta^{\pm})^2$,
 the $(X,Y)$'s in (ii) satisfy (\ref{eq eins sec 5 a}) and (\ref{eq eins sec 5 b}).
 It is a simple check that the discriminant of the quadratic equation is positive.
\end{proof}

\begin{rmk}
We shall refer to the above zeros as the Einstein points $(X^{\pm}_E,Y^{\pm}_E)$, where
the $\pm$ signs correspond to those occurring in $\beta^{\pm}$. Note also that in
\cite{Wa92} only one Einstein metric was found when $m \geq 3$.
\end{rmk}

The following estimates are needed for the discussion of eigenvalues below.
\begin{lemma}
$\beta^{\pm}$ has the following lower  and upper bounds:
 \begin{align}
  & \frac{\sqrt{d}(d+2)}{2 (2d+ 3)} < \beta^{+} < \frac{\sqrt{d}(d+2)}{2 (d+ 1)},
    \label{eq beta pos lower bdd} \\
  &  \frac{\sqrt{d}}{4 (d+ 2)} < \beta^{-} <   \frac{5\sqrt{d}}{6(d+2)}.  \label{eq beta neg  bdd}
 \end{align}
\end{lemma}

\begin{proof}
For the lower bound of $\beta^+$ we compute
\[
\beta^{+} > \frac{d+2}{2\sqrt{d} (2+ 3(1-m^{-1})d^{-1})} > \frac{d+2}{2\sqrt{d} (2+ 3d^{-1})}.
\]
For the upper bound of $\beta^+$ we compute
\[
\beta^{+} < \frac{2(d+2)}{2 \sqrt{d} \left (2 +3(1- m^{-1}) d^{-1} \right )}
\leq \frac{2(d+2)}{2 \sqrt{d} \left (2 +2 d^{-1} \right )},
\]
where we used $m \geq 3$ to get the last inequality.
This proves (\ref{eq beta pos lower bdd}).

For the lower bound of $\beta^-$ we compute using $\sqrt{1- \alpha} \leq 1- \frac{1}{2} \alpha$ for
 $\alpha \in [0,1]$
\begin{align*}
 \beta^{-} & >   \frac{(d+2) \left ( (1+2m^{-1}) d (d+2)^{-2} + \frac{3}{2}(1-m^{-1})(1+2 m^{-1})
(d+2)^{-2} \right ) }{2 \sqrt{d} (2+ 3(1-m^{-1})d^{-1})}  \\
& > \frac{ (d+2)^{-1} \left ((1 +2m^{-1})d + \frac{3}{2} \right) }{2
\sqrt{d} (2+ 3(1-m^{-1})d^{-1})} \\
& > \frac{\sqrt{d}}{4(d+2) },
 \end{align*}
 where we have used $m \geq 3$ to get the last inequality.
For the upper bound of $\beta^-$ we compute using $\sqrt{1- \alpha} \geq 1-  \alpha$ for
 $\alpha \in [0,1]$
\begin{align}
 \beta^{-} & <  \frac{ (d+2)^{-1} \left (2(1+2m^{-1}) d + 3(1-m^{-1})(1+2 m^{-1})  \right )
  }{ 2 \sqrt{d} (2+ 3(1-m^{-1})d^{-1}) }   \notag  \\
 & = \frac{ (d+2)^{-1} }{2 \sqrt{d} } \cdot (1+2m^{-1}) d. \label{eq beta neg  bdd immed}
 \end{align}
The lemma is proved.
\end{proof}

Let  $I_{m \times m}$ be the identity matrix and let $J$ be the $m \times m$ matrix
all of whose entries are $1$.
Define $\mathcal{L}^{\pm}$ to be the linearization of the vector field from system
 (\ref{eq sec 5 ode a}) and (\ref{eq sec 5 ode b}) at the Einstein points $(X^{\pm}_E,Y^{\pm}_E)$.
By a straightforward computation we have
\begin{align}
\mathcal{L}^{\pm} =  \left( \begin{array}{cc}
c^{\pm}_{11}I - 8 (\beta^{\pm})^2 J & c^{\pm}_{12}I - \frac{8\sqrt{d}}{m^2} \beta^{\pm} J  \\
c^{\pm}_{21} I+ \frac{8m^2}{\sqrt{d}} (\beta^{\pm})^3 J  &  c^{\pm}_{22}I + 8 (\beta^{\pm})^2 J
\end{array} \right)_{2m \times 2m} ,
\end{align}
where
\begin{align*}
& c^{\pm}_{11} = 4(2m-1) (\beta^{\pm})^2 -\frac{m}{2},
 &&  c^{\pm}_{12} = \frac{8 \sqrt{d}}{m} \beta^{\pm} >0,  \\
&  c^{\pm}_{21} = \frac{6m^2(m-2)}{ (\sqrt{d})^{3} }\, (\beta^{\pm})^3 >0,
&& c^{\pm}_{22} =  \frac{12(m-1)}{d} \,(\beta^{\pm})^2 - \frac{2(d+2)m}{ \sqrt{d}}\, \beta^{\pm}.
\end{align*}
Note that by the equation for $\beta^{\pm}$ we can rewrite $c^{\pm}_{22} =
\left (\frac{6(m-1)}{d} -4m \right )(\beta^{\pm})^2 -\frac{m+2}{2}$.

We define two matrices
\[
C^{\pm} \doteqdot \left ( \begin{array}{cc}
c^{\pm}_{11} & c^{\pm}_{12} \\ c^{\pm}_{21} & c^{\pm}_{22}
\end{array} \right ) .
\]
\begin{lemma} \label{spectrumC}
The eigenvalues of $C^{\pm}$ are real and distinct.
If we denote by $\lambda_1^{\pm}$ and $\lambda_2^{\pm}$ the two eigenvalues of $C^{\pm}$
in ascending order, then

\vskip .1cm
\noindent $($i$)$  $\lambda_1^+ <0$ and $\lambda_2^+ >0$, and

\vskip .1cm
\noindent $($ii$)$    $\lambda_1^- < \lambda_2^- <0$.
\end{lemma}

\begin{proof} The simplest way to deduce the first statement is to observe that
$C^{\pm}$ can be written as the sum of a positive matrix  and the scalar matrix
$-\frac{2(d+2)}{\sqrt{d}} m \beta^{\pm} I_{2 \times 2}$. The conclusion then follows
from the Perron-Frobenius theorem for positive matrices.

$($i$)$ It suffices to prove that $c_{11}^{+} >0$ and $c_{22}^+ <0$ which then imply
that the determinant $\det C^+ <0$. We consider two subcases.

\noindent Subcase $($i1$)$ $d \geq 2$.  Using (\ref{eq beta pos lower bdd}) we have
\begin{align*}
& c_{11}^{+} >   4(2m-1) \frac{d(d+2)^2}{4 (2d+ 4)^2} -\frac{m}{2} = \frac{1}{4}(2md -d-2m)
>0,  \\
&  c_{22}^{+} < \frac{2(d+2) \beta^+}{\sqrt{d}} \left ( \frac{3(m-1)}{d+1} -m \right ) <0.
\end{align*}

\noindent  Subcase $($i2$)$ $d=1$. By Lemma \ref{lem 5.1 called} we have
$\beta^+ =  \frac{3m+ \sqrt{4m^2 -7m+6}}{2(5m-3)}$. Hence
\[
\frac{5m-2}{2(5m-3)}  < \beta^+ < \frac{5m-1}{2(5m-3)}.
\]
Using these inequalities we can check that $ c_{11}^{+} > 0$ and $ c_{22}^{+} <0$ when $d=1$.

\vskip .1cm
$($ii$)$ First we prove  trace $\tr C^- <0$ by considering two subcases.

\noindent  Subcase $($ii1$)$ $d \geq 2$. We compute  using (\ref{eq beta neg  bdd})
\begin{align*}
\tr C^- & =  \left ( \frac{6(m-1)}{d} +4m-4 \right ) (\beta^-)^2 -(m+1)  \\
& <  \left ( \frac{6(m-1)}{d} +4m-4 \right ) \cdot \frac{25d}{36(d+2)^2} -(m+1).
\end{align*}
The last expression as a function of $d \in [2, \infty)$ achieves its maximum when $d =2$,
and the maximum value is negative. This proves $\tr C^-  <0$ for $d \geq 2$.

\noindent  Subcase $($ii2$)$ $d=1$.  By Lemma \ref{lem 5.1 called} we have
\[
\frac{m+1}{2(5m-3)}  < \beta^- < \frac{m+2}{2(5m-3)}.
\]
From this we can calculate that $c_{11}^-=4(2m-1) (\beta^-)^2 - \frac{m}{2} <0$
and $c_{22}^-= (2m-6)(\beta^-)^2 - \frac{m+2}{2} <0$.
Hence $\tr C^- <0$ when $d=1$.

\medskip

Next we show that $\det C^- >0$ for $m \geq 3$ and any $d$. We compute
\begin{align*}
\det C^- &= \frac{2m\beta^-}{d} \left | \begin{array}{cc}
4(2m-1) (\beta^-)^2 - \frac{m}{2}  &    \frac{8 \sqrt{d}}{m}\beta^-      \\
\frac{3m(m-2)}{\sqrt{d}} (\beta^-)^2      &  \frac{6(m-1)}{m}\beta^-  -(d+2)\sqrt{d}
\end{array} \right |  \\
&= \frac{2m\beta^-}{d} \left (\begin{array}{c}
24(m-1+m^{-1}) (\beta^-)^3 -4\sqrt{d}(d+2)(2m-1) (\beta^-)^2  \\
-3(m-1)\beta^- +\frac{1}{2} \sqrt{d}(d+2) m
\end{array} \right )  \\
& =  \frac{2m\beta^-}{d} \left (\begin{array}{c}
- (8(2m-1)d-12-12m^{-1} ) (\beta^-)^3 \\
- (3(m-1)+(2m+3-2m^{-1})d  )\beta^-  +\frac{1}{2} \sqrt{d}(d+2) m
\end{array} \right )
\end{align*}
where we have used the quadratic equation of $\beta^-$ to get  the last inequality.
We define
\[
\Psi(m,d) \doteqdot   (8(2m-1)d-12-12m^{-1} ) (\beta^-)^2
+ 3(m-1)+(2m+3-2m^{-1})d  .
\]

On the other hand  when $m \geq 3$ and $d \geq 5$ we have
\[
\left (3d^2+2d -\frac{527}{9} \right )m >15d,
\]
which implies
\[
 \left ( \frac{127}{9} m  +(2m+3)d \right ) \cdot \frac{5\sqrt{d}}{6(d+2)}
 <  \frac{1}{2} \sqrt{d}(d+2) m.
\]
By the upper bound of $\beta^-$ we get
\[
 \left ( \frac{127}{9} m  +(2m+3)d \right ) \beta^-
 <  \frac{1}{2} \sqrt{d}(d+2) m.
\]
Combining this with the fact that for $d \geq 2$ and $m \geq 3$, we have
\begin{align*}
\Psi(m,d) < & \, 16md \cdot \frac{25d}{36(d+2)^2} + 3m+(2m+3)d  \\
 <& \,\frac{127}{9}\, m +(2m+3)d,
\end{align*}
and so we deduce that  $-\Psi(m,d) \beta^-  +\frac{1}{2} \sqrt{d}(d+2) m >0$.
Hence $\det C^-  >0$  when $m \geq 3$ and $d \geq 5$.

We have also shown that  $\det C^- >0$ for $m \geq 3$ and $d =1,2,3,4$.
This is done by using the specific value for $d$ and analysing the resulting
expression for $\det C^{-}$. The calculations are straightforward but lengthy
and not very illuminating so we omit them.

From $\tr C^- <0$ and $\det C^- >0$ we conclude that the eigenvalues of $C^-$ are
both negative.
\end{proof}

Note that $\mathcal{L}^{\pm}$ maps linear subspace $V_1 \doteqdot
\{(X,Y) \in \mathbb{R}^{2m},  \sum_{k=1}^m X_k =\sum_{k=1}^m Y_k =0 \}$ into itself.
Using the fact that $c^{\pm}_{12} -  \frac{8 \sqrt{d}}{m} \beta^{\pm} =0$ it is easy to check that
$\mathcal{L}^{\pm}$ in addition maps the linear subspace $V_2 \doteqdot
\{(X,Y) \in \mathbb{R}^{2m},  \sum_{k=1}^m X_k=0 \}$ into itself.
This last property also follows from the fact that the affine hypersurface $\sum_{k=1}^m X_k =1$ is
preserved by (\ref{eq sec 5 ode a}).
For the remainder of this section, we will regard $\mathcal{L}^{\pm}$ as a linear
endomorphism of $V_2$ without introducing a separate notation for its restriction
to $V_2$.

\begin{lemma} \label{lem Linearized L pm eigenvalues}
$($i$)$  The matrices  $\mathcal{L}^{\pm}$ are diagonalizable. They have distinct eigenvalues
 $\lambda_1^{\pm}$ and $\lambda_2^{\pm}$  and the corresponding eigenvectors
form a basis of the linear subspace $V_1$.  $\mathcal{L}^{\pm}$ have a further
eigenvalue $\lambda^{\pm}_3\doteqdot
 \frac{2(d+2)m}{\sqrt{d}}  \beta^{\pm} -(m+2)$ with associated eigenvector lying in $V_2$
but not in $V_1$. Furthermore, $\lambda^{+}_3 >0$ while $\lambda^{-}_3 <0$.

\vskip .1cm
\noindent $($ii$)$ $\mathcal{L}^{+}$ has $m$  positive eigenvalues and $m-1$
    negative eigenvalues on $V_2$.

\vskip .1cm
\noindent  $($iii$)$ $\mathcal{L}^{-}$ is negative definite on $V_2$.
\end{lemma}

\begin{proof}
$($i$)$ Let $\vec{u} =(u_1, \cdots, u_m)^T$ be any column vector which satisfies
$\sum_{i=1}^m u_i =0$. Then we have

\begin{align*}
\mathcal{L}^{\pm}  \left ( \begin{array}{c}
\sigma_1 \vec{u} \\    \sigma_2 \vec{u}
\end{array} \right )
= \left ( \begin{array}{c} (c^{\pm}_{11} \sigma_1 + c^{\pm}_{12}\sigma_2 ) \vec{u}  \\
(c^{\pm}_{21} \sigma_1 + c^{\pm}_{22}\sigma_2 ) \vec{u}
\end{array} \right )
=\lambda   \left ( \begin{array}{c}
\sigma_1 \vec{u} \\ \sigma_2 \vec{u}
\end{array} \right )
\end{align*}
provided that the column vector $( \sigma_1,\sigma_2 )^T \in \mathbb{R}^2$ is
an eigenvector of matrix $C^{\pm}$ with eigenvalue $\lambda$:
\[
C^{\pm}  \left ( \begin{array}{c}
\sigma_1  \\ \sigma_2
\end{array} \right )
   =  \lambda \left ( \begin{array}{c}
\sigma_1 \\ \sigma_2
\end{array} \right ).
\]
Hence by Lemma \ref{spectrumC} we obtain distinct eigenvalues
$\lambda_1^{\pm}$ and $\lambda_2^{\pm}$ of $\mathcal{L}^{\pm}$
and  the associated eigenvectors form a basis of the linear subspace $V_1$.

 Next, let $\vec{e}$ be the column vector $(1,\cdots,1 )^T$. We compute that
\begin{align*}
\mathcal{L}^{\pm}  \left ( \begin{array}{c}
\sigma_1 \vec{e} \\    \sigma_2 \vec{e}
\end{array} \right )
= \left ( \begin{array}{c} (c^{\pm}_{11} \sigma_1 - 8m (\beta^{\pm})^2 \sigma_1
+ c^{\pm}_{12}\sigma_2 - \frac{8 \sqrt{d}}{m}\beta^{\pm} \sigma_2 ) \vec{e}  \\
(c^{\pm}_{21} \sigma_1+\frac{8m^3}{\sqrt{d}} (\beta^{\pm})^3 \sigma_1 +
 c^{\pm}_{22}\sigma_2 +8m (\beta^{\pm})^2 \sigma_2 ) \vec{e}
\end{array} \right )
=\lambda   \left ( \begin{array}{c}
\sigma_1 \vec{e} \\ \sigma_2 \vec{e}
\end{array} \right )
\end{align*}
provided that the  column vector $( \sigma_1,\sigma_2 )^T \in \mathbb{R}^2$
is an eigenvector with eigenvalue $\lambda$ of matrix $D^{\pm}$:
\[
D^{\pm}  \left ( \begin{array}{c}
\sigma_1  \\ \sigma_2
\end{array} \right )
 \doteqdot \left ( \begin{array}{cc}
c^{\pm}_{11} - 8m (\beta^{\pm})^2  & c^{\pm}_{12} - \frac{8 \sqrt{d}}{m}\beta^{\pm}
\\ c^{\pm}_{21} +\frac{8m^3}{\sqrt{d}} (\beta^{\pm})^3  & c^{\pm}_{22} + 8m (\beta^{\pm})^2
\end{array} \right )
 \left ( \begin{array}{c}
\sigma_1  \\ \sigma_2
\end{array} \right )   =  \lambda \left ( \begin{array}{c}
\sigma_1 \\ \sigma_2
\end{array} \right ).
\]
Since
\[
D^{\pm}= \left ( \begin{array}{cc}
- 4 (\beta^{\pm})^2 -\frac{m}{2}  & 0  \\
\frac{6m^2(m-2) +8m^3d}{(\sqrt{d})^{3}} \,(\beta^{\pm})^3  & \frac{ 8md +12(m-1)}{d}\, (\beta^{\pm})^2
-\frac{2(d+2)m}{\sqrt{d}}\, \beta^{\pm}
\end{array} \right ) ,
\]
one sees that the vector $( \sigma_1,\sigma_2 )^T =(0,1)^T$ is an eigenvector with eigenvalue
\begin{equation}  \label{lambda3}
  \lambda_3^{\pm} =
\frac{ 8md +12(m-1)}{d}\, (\beta^{\pm})^2  -\frac{2(d+2)m}{\sqrt{d}}\, \beta^{\pm} =
\frac{2(d+2)m}{\sqrt{d}}\, \beta^{\pm} -(m+2).
\end{equation}
This gives an eigenvector in the subspace $V_2$ which does not lie in $V_1$.

Note that using (\ref{eq beta pos lower bdd}) we get
\[
\lambda^{+}_3 > \frac{(d+2)^2}{2d+3}m-(m+2) >0 ,
\]
while using (\ref{lambda3}) and (\ref{eq beta neg  bdd immed})
yields $\lambda^{-}_3 < 0$.

\vskip .1cm
$($ii$)$ and  $($iii$)$ These follow  immediately from $($i$)$.
\end{proof}

 \subsection{Ancient solutions on $\bar{P}$ with $m\geq 3$}
Now we can prove the existence of ancient solutions on the special class of fibre bundles
$\bar{P}$ with $m \geq 3$ which were described at the beginning of \S 5.

\begin{theorem} \label{thm sec 5 m> 3}
Let  $(M_i^{4n_i}, g_i), i=1, \cdots,m, $ be QK manifolds with $\operatorname{Rc}_{g_i} =\Lambda_i g_i$.
 Assume $m \geq 3$, $n_1 =\cdots = n_m =d$, and $\Lambda_i =d+2$ for all $i$.
 Let  $g^+$ and $g^-$ be the Einstein metrics given in Lemma \ref{lem 5.1 called}.

 \vskip .1cm
\noindent  $($i$)$ There is a continuous $(m-2)$-parameter family of ancient solutions
of the Ricci flow on fibre bundle $\bar{P}$ such that
$\lim_{\tau \rightarrow \infty} \frac{1}{\tau} \bar{g}(\tau)$ is a
multiple of the Einstein metric $g^+$.

 \vskip .1cm
\noindent  $($ii$)$ There is a continuous $(2m-2)$-parameter family of ancient solution of
the Ricci flow such that $\lim_{\tau \rightarrow \infty}  \frac{1}{\tau} \bar{g}(\tau)$ is a
multiple of the Einstein metric $ g^-$.

 \vskip .1cm
\noindent  $($iii$)$ All these ancient solutions are $\kappa$-noncollapsing, of type I as
 $\tau \rightarrow \infty$, and have positive Ricci curvature when $\tau$ is large enough.
\end{theorem}

\begin{proof}
$($i$)$ We consider solutions $(X,Y)$ of the $\operatorname{ODE}$ system
(\ref{eq sec 5 ode a}) and (\ref{eq sec 5 ode b}) in the space $\{(X,Y) \in \R^{2m}: \sum_{k=1}^m X_k=1 \}$.
By Lemma \ref{lem Linearized L pm eigenvalues} the linearized operator $\mathcal{L}^+$ at the
Einstein point $(X^+_E,Y^{+}_E)$ has a negative eigenvalue of multiplicity $m-1$ in $V_2$.
So we may apply the Hartman-Grobman Theorem to conclude that there is a  continuous $(m-2)$-parameter
family of  solutions $(X^+(u),Y^+(u))$ whose limits $\lim_{u \rightarrow \infty}(X^+(u),Y^+(u)) =
 (X^+_E,Y^+_E)$.
 By the discussion after equation (\ref{eq sec 5 ode b}) we get  a  continuous
  $(m-2)$-parameter family of  solutions of the backwards Ricci flow
   (modulo time translation and parabolic scaling).
 To see that these solutions are ancient, we prove next that
  $\lim_{u \rightarrow \infty}\tau(u) \rightarrow \infty$.

If we take a flow line $(X^+(u),Y^+(u))$ in the local stable manifold near $(X^+_E,Y^+_E)$,
then  $\lim_{u \rightarrow \infty} (X^+(u), Y^+(u)) =(X^+_E,Y^+_E)$
and $\lim_{u \rightarrow \infty} E(X^+(u),Y^+(u)) = E(X^+_E, Y^+_E)$.
Hence given any $\varepsilon >0$ there is a $u_0$ such that for $u \geq u_0$ we have
\[
 |E(X^+(u), Y^+(u)) - E(X^+_E, Y^+_E)| \leq  \varepsilon  .
\]
We compute by using  $\frac{d}{du} \ln \hat{a}(u) =E(X^+(u), Y^+(u))$ that
\[
 \ln \hat{a}(u) -\ln \hat{a}(u_0) \leq ( E(X^+_E, Y^+_E)
+ \varepsilon) (u -u_0) \quad \text{ for }  u \geq u_0.
\]
Hence we have
\[
 \hat{a}(u) \leq \hat{a}(u_0) e^{ ( E(X^+_E, Y^+_E) +\varepsilon) (u -u_0)}
 \quad \text{ for } u \geq u_0.
\]
This and the relation $d \tau = \hat{a}(u) du$
imply that $\lim_{u \rightarrow \infty}\tau(u) \rightarrow \infty$.

From  (\ref{eq quaternion a evol}) and (\ref{eq quaternion b})  we get
\begin{align*}
& \frac{d a_k}{d \tau} = \frac{1}{2} +X_k +4dX_k^2Y_k^2 \rightarrow \frac{m+2}{2m}
+ 4(\beta^{+})^2   \\
& \frac{d b_k}{d \tau} = 2(d+2)-6X_k(1-X_k)Y_k   \rightarrow 2(d+2) -\frac{6(m-1)}{m\sqrt{d}}
\beta^{+}  >0
\end{align*}
as $\tau \rightarrow \infty$. We have used (\ref{eq beta pos lower bdd})
to get the last inequality. This implies that $\lim_{\tau \rightarrow \infty} \frac{a_k}{\tau}$ and
 $\lim_{\tau \rightarrow \infty} \frac{b_k}{\tau}$ exist. $($i$)$ is now proved.

\noindent{$($ii$)$} Note that  the linearized operator $\mathcal{L}^-$ at the
Einstein point $(X^-_E,Y^{-}_E)$
has $2m-1$ negative eigenvalues in $V_2$ by Lemma \ref{lem Linearized L pm eigenvalues}.
With this observation, the rest of the proof is similar to the proof of $($i$)$, and so
we omit it.

\noindent  $($iii$)$ The $\kappa$-noncollapsing and Type I properties follow from the fact that
 $\lim_{\tau \rightarrow \infty} \frac{\bar{g}(\tau)}{\tau}$ is a multiple of
 the nonflat Einstein metrics $g^{\pm}$ and Corollary \ref {cor used to get type I asypm}.

  From \cite[(2.2)]{Wa92}, the restriction of the Ricci tensor of $\bar{g}(\tau)$
to the fibre directions is positive definite. So by (\ref{eq quaternion b})) the
Ricci tensor is positive definite if and only if
\[
 2(d+2) -6X_k(1-X_k)Y_k>0, \quad \text{for } k=1, \cdots, m.
\]
These inequalities hold when $(X(\tau),Y(\tau))$  are close to $(X^{\pm}_E,Y^{\pm}_E)$.
\end{proof}

\begin{rmk}
$($i$)$ The dimensions of the bundles considered in this section are of the form
$4md+3(m-1)$. In particular, for $m=2l+1$ with $ l \geq 1$, the dimension is even.
However, the total spaces of the bundles cannot admit a K\"ahler structure if none
of the QK factors in the base are complex Grassmannian. To see this, recall that
S. Salamon \cite[14.83]{Be87} has proved that positive QK manifolds are simply connected,
and that the rigidity result of LeBrun and Salamon \cite{LS94} implies that under our
assumption the second Betti number of the base is $0$. As well, the fibres are
products of $\R\PP^3$, which is a homology $3$-sphere. It then follows for example
 from the Serre spectral sequence that the second Betti number of the total space
 is trivial. So the metrics in our ancient flows do not admit compatible symplectic
structures which are non-cohomologous to zero. The only other ancient solutions
of this type that we know of are the continuous families on compact simple Lie groups
constructed by Lauret \cite{La13} and the isolated solutions on the homogeneous spaces
$\Sp(3)/(\Sp(1) \times \Sp(1) \times \Sp(1))$ and ${\rm F}_4/ \Spin (8)$
(see Example 5 in \cite{BKN12}).

$($ii$)$ We also obtain  the pseudo-Riemannian  Ricci flow analogs
of the ancient solutions in Theorem \ref{thm sec 5 m> 3} by the same arguments
in Remark \ref{pseudo-Rflow}.
\end{rmk}


\end{document}